%% file: totally_isotropic_arxiv_v2.1.tex
\documentclass[11pt,a4paper, notitlepage]{article}
\usepackage{import} %This is for importing
\import{}{preamble}

\makeatletter                                                                                                                             

\theoremstyle{plain}
\newtheorem{pl}[thm]{Presentation Lemma~}

\usepackage{footnote}
\newcommand{\ccite}[2]{\cite[#2]{#1}}
\newcommand{\stacks}[1]{\ccite{stacks-project}{\href{https://stacks.math.columbia.edu/tag/#1}{Tag \MakeUppercase{#1}}}}

\title{\uppercase{\textbf{\large{Isotropic Torsors on Smooth Algebras over Prüfer Rings}}}}
\author{\uppercase{Arnab Kundu}{\let\thefootnote\relax\footnote{University of Toronto, St.~George Campus, Toronto, Canada\newline\hspace*{0.48cm} Email: arnab.kundu@utoronto.ca.}}}
\date{}
\usepackage{doi}
\begin{document}
\maketitle
\abstract{The Grothendieck--Serre conjecture predicts that every generically trivial torsor under a reductive group over a regular semilocal ring is itself trivial. Extending the work of Česnavičius and Fedorov, we prove a non-noetherian analogue of this conjecture for rings $A$ that are semilocalisations of smooth schemes over valuation rings of rank one, and for reductive $A$-group schemes $G$ that are totally isotropic. Roughly speaking, such group schemes are characterised by the existence of a parabolic subgroup of their adjoint quotients. Since quasi-split groups are totally isotropic, our result, in particular, generalises the Grothendieck--Serre result of Guo--Liu and the author's thesis. Our proof relies on a new instance of Gabber's presentation lemma, obtained by extending techniques developed in the author's thesis.}
\tableofcontents
{\let\thefootnote\relax\footnote{June 2025}}
%%%%%%%%%%%%%%%%%%%%%%%%%%%%%%%%%%%%%%%%%%%%%%%%%%%%%%%%%%%%%%%%%%%%%%%%%%%%%%%%%%%%%%%%%%%%%%%%%%%%%%%%%
%	We prove descent results in various papers from \cite{antieau-mathew-morrow} and \cite{gillet_levine}
\section{An analogue of the Grothendieck--Serre conjecture for non-noetherian rings.}
	One of the central problems in the study of torsors under reductive group schemes is the Grothendieck--Serre conjecture \ccite{ces_grothendieck-serre}{Conjecture~1.1}, originating from the Chevalley seminar papers of Serre in \ccite{serre_chevalley_58}{page 31, remarque} and Grothendieck in \ccite{grothendieck_chevalley_58}{pages 26–27, remarques 3}.
	\par This conjecture, which may be viewed as a non-abelian analogue of Gersten's injectivity conjecture for algebraic $K$-theory, asserts that any generically trivial torsor over a regular semilocal ring is trivial. As such, it bridges geometry with arithmetic by predicting that $G$-torsors, inherently geometric objects, may be represented by classes in Galois cohomology, a fundamentally arithmetic invariant. %\rl{This conjecture has wide applications blah blah blah. WE can think of this conjecture as a an abelian variant of the Gersten's conjecture for K-theory.} 
	\par The purpose of this article is to investigate the following variant of the Grothendieck--Serre conjecture.
	\bcj[Česnavičius, \ccite{phd-thesis}{Conjecture~2.1.2}]\label{conj:gro_serre_valuation}
		Let $V$ be a valuation ring, let $A$ be an integral domain that is the semilocalisation of a smooth $V$-algebra at finitely many points and let $K$ be the fraction field of $A$. Given a reductive $A$-group scheme $G$, any generically trivial $G$-torsor $E$ over $A$ is trivial, i.e., the restriction morphism induces an injection \bud \ker(H^1(A,G)\to H^1(K, G))=\{\ast\}.\eud
	\ecj
	\par We are particularly interested in the conjecture above for non-noetherian valuation rings $V$. However, before proceeding with our discussion, we note that when $V$ is instead noetherian, Conjecture~\ref{conj:gro_serre_valuation} already recovers two important cases of the Grothendieck--Serre conjecture. 
	\br If $V$ is a field, Conjecture~\ref{conj:gro_serre_valuation} is the equicharacteristic case of the Grothendieck--Serre conjecture, proven by \citeauthor{fedorov_panin} in \cite{fedorov_panin}. Meanwhile, for mixed-characteristic discrete valuation rings $V$, it corresponds to the unramified case, studied by \citeauthor{ces-fedorov} in \cite{fedorov_gro-serre_mixed_published22, fedorov_gro-serre_unramified_mixed22, ces_grothendieck-serre} and \cite{ces-fedorov}. The current state-of-the-art is that this mixed-characteristic and unramified case is known for totally isotropic $G$, i.e., a reductive $A$-group scheme $G$ such that its adjoint quotient $G^{\text{ad}}$ has no anisotropic factors (see Definition~\ref{defn:totally_istropic_groups}).\er %In mixed-characteristic, the cases beyond the unramified one currently remains out of reach.
	\par The expectation that Conjecture~\ref{conj:gro_serre_valuation} holds true stems from the fact that smooth algebras over valuation rings behave as non-noetherian analogues of regular rings. This hypothesis is reinforced by Zariski's local uniformisation conjecture (see, for example, \ccite{phd-thesis}{Conjecture~2.1.1}), which predicts that each smooth algebra over a valuation ring is ind-regular local. In particular, by a limit argument, Conjecture~\ref{conj:gro_serre_valuation} follows from the Grothendieck--Serre conjecture combined with Zariski's conjecture. However, Zariski's conjecture--which is implied by the resolution of singularities--remains vastly open in positive and mixed-characteristic.
	\par In this article, we demonstrate a case of Conjecture~\ref{conj:gro_serre_valuation} for valuation rings $V$ of rank one without employing Zariski's local uniformisation conjecture. More precisely, \textsection\ref{section:main_proof} is dedicated to the proof of the following. 
	\par Recall that Prüfer rings (see Definition~\ref{defn:prufer_domain}) are non-noetherian analogues of Dedekind rings.%, and are characterised by the fact that each of their local rings is a valuation ring. However, unlike Dedekind domains, they are allowed to be of arbitrary Krull dimensions.
	\bt\label{thm:GS-totally_isotropic}
		Let $R$ be a Prüfer ring of Krull dimension one and $A$ be the semilocalisation of a smooth $R$-algebra at finitely many points. Then, for any totally isotropic reductive $A$-group scheme $G$, every generically trivial $G$-torsor $E$ over $A$ is trivial.%, i.e., we have \bud \ker(H^1(A,G)\to H^1(K, G))=\{\ast\}.\eud
	\et
	 Prior to this work, the case of Theorem~\ref{thm:GS-totally_isotropic} in which $G$ contains a Borel $R$-subgroup was established by \citeauthor{ning_liu_quasi-split} in \cite{ning_liu_quasi-split}. A related subcase--where, in addition, $R$ is a valuation ring of rank one--was obtained simultaneously and independently in the author's thesis in \cite{phd-thesis}.
	\par Our proof (see \textsection\ref{intro:sketch} for a sketch) of Theorem~\ref{thm:GS-totally_isotropic} in \textsection\ref{section:main_proof} builds on the techniques of \cite{phd-thesis} and \cite{ces-fedorov}, and closely follows the strategy developed by \citeauthor{ces-fedorov} in op.~cit. The key geometric input is established in \textsection\ref{section:presentation_lemma}; it is a version of Gabber’s presentation lemma over Prüfer domains, as developed in \cite{phd-thesis} and \cite{amar:unramified_case_gersten_conjecture} and recalled below.
	\begin{pl}\label{prop:presentation_lemma}
		For 
		\bun
		\item[$\circ$] a Pr\"ufer ring $R$ of Krull dimension $\leqslant 1$, 
		\item[$\circ$] a smooth $R$-scheme $X$ fibrewise of pure relative dimension $d>0$,
		\item[$\circ$] a closed subscheme $Z\subset X$ that is of codimension $\geqslant 2$,\footnote{This is equivalent to the condition that the $R$-generic fibres of $Z$ are of codimension $\geqslant 2$ while its other $R$-fibres are of codimension $\geqslant 1$ in that of $X$.}
		\item[$\circ$] a closed subscheme $Z\subset Y\subset X$ that is $R$-fibrewise of codimension $\geqslant 1$, and
		\item[$\circ$] points $x_1,\ldots,x_n\in X$;
		\eun 
		there are affine opens $x_1,\ldots,x_n\in U\subset X$ and $S\subset \bb{A}^{d-1}_{R}$ and a smooth $R$-morphism $\pi\colon U\to S$ of pure relative dimension $1$ such that $\pi|_{Z\cap U}$ is finite and $\pi|_{Y\cap U}$ is quasi-finite. 
	\end{pl}
	 We note two key differences between Gabber's presentation lemma (see \cite{ct_hoobler_kahn} and \cite{hogadi_kulkarni_gabber_finite_field}) and Presentation~Lemma~\ref{prop:presentation_lemma}. 
	 \par First, in order to ensure that $Z\cap U$ is $\pi$-finite, our method requires the assumption that $Z\subset X$ is of codimension $\geqslant 2$. This is a restrictive condition that reflects the inherent difficulties of working in mixed characteristic. As a consequence, we are only able to guarantee that $Y\cap U$ is $\pi$-quasi-finite (the case when $Z=\emptyset$ was established in \ccite{amar:unramified_case_gersten_conjecture}{Presentation~Lemma~3.2}). 
	 \par Second, we do not construct a closed embedding $Z\hookrightarrow\bb{A}^1_R$, as it is not necessary for our application. Instead, we rely on ~\ccite{ces-fedorov}{Lemma~2.5}. As shown by \citeauthor{schmidt_strunk} in \cite{schmidt_strunk}, such a closed embedding can be constructed Nisnevich-semilocally around $x_1,\ldots x_n\in X$, at least, when $R$ is a Dedekind domain whose residue fields are all infinite. We plan to pursue their approach in a forthcoming project with T.~Bouis, where we aim to demonstrate a Nisnevich-local presentation lemma over arbitrary Prüfer domains.
	 \par To prove Presentation Lemma~\ref{prop:presentation_lemma} in \textsection\ref{section:presentation_lemma}, we adapt the approach from \ccite{ces_grothendieck-serre}{Proposition~4.1} to our setting. The idea is to slice a compactification $X\hookrightarrow\overline{X}\subset\bb{P}^{n}_R$ by $d-1$ hypersurfaces $H_1,\ldots,H_{d-1}$ in general positions. These hypersurfaces are selected to intersect the boundary $\overline{X}\setminus X$ in a controlled way. Specifically, we ensure that \bud Y\cap H_1\cap\cdots\cap H_{d-1}\text{ is finite and }(\overline{Y}\setminus Y)\cap H_1\cap\cdots\cap H_{d-1}=\emptyset\eud The rational morphism $\pi\colon X\dashrightarrow\bb{A}^{d-1}_V$ is then defined by $H_1,\ldots, H_{d-1}$. By ensuring that $\pi$ is smooth at each $x_1,\ldots,x_n$, we find a neighbourhood  $x_1,\ldots,x_n\in U\subset X$ where $\pi$ is smooth. The $\pi$-finiteness of $Z\cap U$ then follows by the properness of $\overline{Z}$ and the $\pi$-quasi-finiteness of $Z\cap U$.
	 \par The case where $R$ has infinite residue fields is simpler, as we can arrange for $H_1,\ldots,H_d$ to be hyperplanes. In this situation, the hyperplane sections obtained by using Bertini's theorem at the $R$-special fibres of $X$ can be lifted globally. However, when $R$ has a finite residue field, ensuring that $H_1,\ldots,H_{d-1}$ have the same degree becomes challenging. This is due to the nature of Bertini's theorem over finite fields (see \cite{poonen_bertini_finite_field}, cf.~\ccite{ces_grothendieck-serre}{Lemma~3.2}). As a consequence, these hypersurfaces $H_1,\ldots,H_{d-1}$, which are possibly of different degrees, determine a rational morphism to a weighted projective space. This introduces additional complexity, necessitating the use of a weighted blowup to lift the resulting rational morphism. For more details on weighted projective spaces and weighted blowups, we refer the reader to \ccite{ces_grothendieck-serre}{\textsection6.1}.
	 %\par The primary challenge in proving \Cref{lem:intro_presentation_lemma} arises from the intricate algebraic structure of such rings $R$, which are non-noetherian, making their algebraic geometry complex. In \textsection\ref{section:presentation_lemma}, we develop the necessary background on the algebraic geometry of Prüfer rings. When $R$ has Krull dimension 1, we can still manage the behavior of $R$-schemes $X$, as these schemes only have two types of $R$-fibers: special and generic. This is why the Krull dimension assumption on $R$ is required in \Cref{lem:intro_presentation_lemma}, and consequently, in \Cref{thm:GS-totally_isotropic}. 
	 \par
	 \bmo[Sketch of our proof of Theorem~\ref{thm:GS-totally_isotropic}]\label{intro:sketch}~
	\bn[(1)]
	\item We spread out $A$ to an affine $R$-smooth scheme $X$ such that $G$ and $E$ are defined over $X$. By a standard argument, we produce a closed subscheme $Y\subset X$ that is $R$-fibrewise of positive codimension such that $E$ trivialises away from $Y$. Thanks to Presentation~Lemma~\ref{prop:presentation_lemma} (in fact, \ccite{amar:unramified_case_gersten_conjecture}{Presentation~Lemma~3.2} is enough), purity for $G$-torsors \ref{prop:purity_for_torsors} and the generalised Horrocks' principle \ref{lem:generalised_horrocks_principle}, we find a $G$-torsor $\scr{E}_X$ over $\bb{P}^1_{X\setminus Z}$, where $Z\subset X$ is a closed subscheme contained in $Y$ and is of codimension $\geqslant 2$, such that $\scr{E}_X|_{\{t=0\}}\cong E|_{X\setminus Z}$, and at the same time, $\scr{E}_X|_{\{t=\infty\}}$ and $\scr{E}_X|_{\bb{P}^1_{X\setminus Y}}$ are trivial (see Proposition~\ref{prop:2.6}). 
	\item By another application of Presentation~Lemma~\ref{prop:presentation_lemma}, pulling-back the data from $X$, we have an open $C\subset\bb{A}^1_A$ equipped with an $A$-section $s$, as well as an $A$-quasi-finite closed subscheme $\scr{Y}\subset C$ and an $A$-finite closed subscheme $\scr{Z}\subset\scr{Y}$. Furthermore, by \ccite{ces-fedorov}{Proposition~3.1(a)}, there are a $G$-torsor $\scr{E}$ over $C$ that trivialises away from $\scr{Y}$ along with a $G$-torsor $\tilde{\scr{E}}$ over $\bb{P}^1_{C\setminus\scr{Z}}$ such that $\tilde{\scr{E}}|_{\{t=0\}}\cong \scr{E}|_{C\setminus\scr{Z}}$, and at the same time, $\tilde{\scr{E}}|_{\{t=\infty\}}$ and $\tilde{\scr{E}}|_{\bb{P}^1_{C\setminus\scr{Y}}}$ are trivial (see Theorem~\ref{thm:GS-totally_isotropic_body}).
	\item Finally, again by the generalised Horrocks' principle \ref{lem:generalised_horrocks_principle}, we find an $A$-finite closed subscheme $\scr{Z}\subset H\subset C$ such that $\scr{E}$ trivialises away from $H$, and hence, $\scr{E}$ extends all the way to a torsor over $\bb{P}^1_A$. However, by the sectional invariance \ref{lem:section_invariance}, since $\scr{E}|_{\{t=\infty\}}$ is trivial, the same must be true for $E\cong\scr{E}|_{\{t=0\}}$. This concludes the proof.
	%\item Fix PRüfer domain notation as $R$???
	\en
	\emo
%	\bmo{\bfseries Organisation of this article}\label{sketch} \emo
%	Section \ref{section:presentation_lemma} is dedicated to the proof of Presentation Lemma~\ref{prop:presentation_lemma} and Section \ref{section:main_proof} is dedicated to the same of \Cref{thm:GS-totally_isotropic}.
	\subsection*{\bfseries Notations and conventions}
		Let $S$ be a scheme, $s\in s$ be a point and $f\colon S'\to S$ be a morphism of schemes.  
		\par
		\bun			
		\item[$\circ$] Then, the localisation (resp., the residue field) of $S$ at $s$ shall be denoted by $\ca{O}_{S,s}$ (resp., by $\kappa(s)$). %$\Sets_{\ast}$ denotes the category of pointed sets.
		\item[$\circ$] When $S=\spec(A)$ is affine, the residue field at a prime ideal $\fr{p}\subset A$ shall be denoted by $\kappa(\fr{p})$. %We shall denote the $I$-adic completion $\lim_nA/I^n$ of $A$ by $\widehat{A_{}}^I$ or abusively by $\widehat{A}$, when $I$ is clear from the context. We say that $A$ is \textit{$I$-adically complete} if it is $I$-adically complete and separated, i.e., if the canonical map is an isomorphism $A\iso \lim_nA/I^n$.
		\item[$\circ$] The base change of an $S$-scheme $X$ along $f$ shall be denoted by $X_{S'}$.%The vanishing locus of $I$ is denoted by $V(I)\subseteq\spec A$. If $I$ is principal with generator, say $t$, then $V(t)$ denotes $V(I)$.
		\item[$\circ$]  The $0$-section (resp., $1$-section, resp., $\infty$-section) of $\bb{P}^1_S$ shall be denoted by $\{t=0\}$ (resp., $\{t=1\}$, resp., $\{t=\infty\}$).%The maximal spectrum of $A$ is denoted by $\mathrm{MaxSpec}(A)\subseteq\spec A$.
		%\item[$\circ$] If $I$ is prime, then the localisation (resp., the residue field) of $A$ at $I$ shall be denoted by $A_I$ (resp., by $\kappa(I)$).
		%\item[$\circ$] If $A$ is an integral domain, then the fraction field of $A$ is denoted by $\Frac(A)$.
		%\item[$\circ$] For an $A$-group scheme $G$, the unipotent radical of $G$ is denoted by $\scr{R}_u(G)$.
		%\item[$\circ$] $A'$ is called \textit{essentially smooth} if it can be obtained as the semilocalisation at finitely many primes of a smooth $A$-algebra.
		%\item[$\circ$] The base change of $X$ along $a$ is denoted by $X_{A'}$.
		%\item[$\circ$] The topological space associated to $X$ is denoted by $|X|$.
		\eun
	\subsection*{\bfseries Acknowledgements}
		%The author extends gratitude to Matthew Morrow and Elden Elmanto for their encouragement, which inspired the author to pursue Gersten's injectivity conjecture in the non-Noetherian case. 
		Special thanks are extended to Tess Bouis, Neeraj Deshmukh, Ofer Gabber, Ritankar Nath and David Rydh for many insightful conversations. This project began during the author’s doctoral studies, and he remains intellectually indebted to Elden Elmanto and K\k{e}stutis \v{C}esnavi\v{c}ius for their unwavering support throughout. 
		\par This project received partial funding from the European Research Council under the European Union’s Horizon 2020 research and innovation program (grant agreement No. 851146) and NSERC Discovery grant RGPIN-2025-07114, ``Motivic cohomology: theory and applications''.
\section{Gabber's presentation lemma over stably coherent rings.}\label{section:presentation_lemma}
	In this section, our goal is to prove Presentation Lemma~\ref{prop:presentation_lemma}, which is a version of Gabber's presentation lemma over Prüfer rings of Krull dimension $\leqslant 1$. This result simultaneously generalises both \ccite{phd-thesis}{Lemma~6.4} as well as \ccite{amar:unramified_case_gersten_conjecture}{Presentation Lemma~3.2}. It will be deduced as a corollary of a stronger presentation lemma--namely, Theorem~\ref{prop:moving_lemma_general}--that applies even over stably coherent rings (see Definition~\ref{defn:coherent_rings}).  
	\par Two auxiliary results play a crucial role in this development. First, in Proposition~\ref{lem:valuation_rings_catenary}, we demonstrate that Pr\"ufer domains of finite Krull dimension are universally catenary (we recall the definition from \stacks{00nl} below). This property plays a key role in the dimension-counting argument underlying the proof of Presentation Lemma~\ref{prop:presentation_lemma}. 
	\par Second, in Lemma~\ref{lem:ega_III_1_2.2.2}, we establish that pushforward of coherent sheaves along projective morphisms remains coherent. It is crucial for lifting sections of line bundles from the special fibre in the proof of Theorem~\ref{prop:moving_lemma_general}. While writing this article, we realised that there is a more general version of this result in \ccite{fujiwara_kato_foundations_rigid_geometry}{Chapter I, Theorem 8.1.3}.
	\par We begin our discussion by recalling the notion of universally catenary.
	\begin{defns}
		A topological space $X$ is called \textit{catenary} if for every pair of irreducible closed subsets $T\subset T'$ there exists a maximal chain of irreducible closed subsets $T=T_0\subset T_1\subset\ldots\subset T_n=T'$ and every such chain has the same length (\cite[\href{https://stacks.math.columbia.edu/tag/02I1}{Tag 02I1}]{stacks-project}). A scheme is called \textit{catenary} if its underlying topological space is catenary (\cite[\href{https://stacks.math.columbia.edu/tag/02IW}{Tag 02IW}]{stacks-project}). A scheme $S$ is called \textit{universally catenary} if any locally of finite type $S$-scheme is catenary. A ring is called \textit{catenary} (resp., \textit{universally catenary}) if its spectrum is catenary (resp., universally catenary).
	\end{defns}
	\par We recall the following minor generalisation of \ccite{egaIV_3}{lemme 14.3.10}, which is used to bound the fibres of finite type schemes over Pr\"ufer domains in the proof of Proposition~\ref{lem:valuation_rings_catenary}.
	\bl\label{lem:valuation_ring_fibres_same_dimension}
	Let $R$ be a Pr\"ufer domain, and let $\gamma,\eta\in\spec R$ be points such that $\eta$ is the generic point. Given an irreducible, finite type, dominant $R$-scheme $X$, if $X_{\kappa(\gamma)}\neq\emptyset$ then $\dim X_{\kappa(\eta)}=\dim X_{\kappa(\gamma)}$.
	\el
	\bs
		Since the statement is local we can localise at $\gamma$ and assume that $R$ is a valuation ring with a closed point $\gamma$. We then apply loc.~cit.~to conclude the proof.
	\es
	\br
		The seemingly surprising claim in Lemma~\ref{lem:valuation_ring_fibres_same_dimension} can be explained by noting that the hypothesis ensures that $X$ is $R$-flat (see \ccite{bourbaki}{Chapter~I, \textsection2.4, Proposition~3(ii)}). This flatness condition plays a crucial role in maintaining the dimension of the fibres across different points.
	\er
	We are now ready to show that Pr\"ufer domains of finite Krull dimension are universally catenary. Our proof will closely follow the argument outlined in \ccite{gabber-ramero-foundations}{Lemma 11.5.8}.
	\bp\label{lem:valuation_rings_catenary}
		Let $R$ be a Pr\"ufer domain of finite Krull dimension and let $f\colon X\to\spec R$ be finite type morphism of schemes. The function \bud \delta\colon |X|\to\bb{Z} \text{, given by } \delta(x)=\trdeg_{\kappa(f(x))}(\kappa(x))-\codim(\overline{\{f(x)\}}),\eud is a `dimension' function (cf.~\cite[\href{https://stacks.math.columbia.edu/tag/02I8}{Tag 02I8}]{stacks-project}), i.e., $x$ specialises to $y\neq x$ only if $\delta(x)>\delta(y)$, and a specialisation $x\rightsquigarrow y$ is immediate if and only if $\delta(x)=\delta(y)+1$. Furthermore, if $Y$ is the spectrum of a semilocalisation of $X$, then $|Y|$ is a catenary topological space of finite Krull dimension.  
	\ep
	%\rl{Need to show that the above quantity $\delta$ is bounded. If yes, then the above as a topological space is Noetherian. Take a infinite descending chain of closed subsets. Since, the space is sober each closed subset has a generic point. Thus taking the sequence of generic points $\{x_1,x_2,\ldots\}$. For each $i$, $x_i$ specialises to $x_{i-1}$. Unless the sequence is eventually constant we get a sequence, $\delta(x_1)>\delta(x_2)>\cdots$ contradicting the fact that it is Noetherian.}
	\bs
		We show that it suffices to assume that $Y=X$ (i.e., it is a semilocalisation of $X$ at the empty set) to prove the final statement. First, we claim that it is enough to show that $X$ is catenary. Indeed, this is true since the semilocalisation of any catenary scheme is catenary (being catenary is a Zariski local property (\cite[\href{https://stacks.math.columbia.edu/tag/02I2}{Tag 02I2}]{stacks-project}) and any localisation of a catenary ring is catenary (\cite[\href{https://stacks.math.columbia.edu/tag/00NJ}{Tag 00NJ}]{stacks-project})). Second, by definition of the Krull dimension, it is enough to check that $|X|$ has finite Krull dimension. Therefore, without loss of generality, we may assume that $Y=X$. Henceforth, we show that $|X|$ is a catenary topological space of finite Krull dimension.
		\par A sober topological space (\cite[\href{https://stacks.math.columbia.edu/tag/004X}{Tag 004X}]{stacks-project}) with a dimension function is catenary (see \cite[\href{https://stacks.math.columbia.edu/tag/02IA}{Tag 02IA}]{stacks-project}). In fact, a sober topological space with a bounded dimension function is of finite Krull dimension. Indeed, consider a descending chain $|X|\supseteq X_0\supsetneq X_1\supsetneq\ldots\supsetneq X_m$ of irreducible closed subsets. For each $n$, let $x_n\in X_n$ be the generic point. The containment $X_n\supsetneq X_{n+1}$ implies that $x_n\rightsquigarrow x_{n+1}$ and $x_n\neq x_{n+1}$, and hence, $\delta(x_n)>\delta(x_{n+1})$. %We choose $x_0$ to be the generic point of an irreducible component $X_0'\subseteq X_0$ such that $X_0'\cap X_1\neq\emptyset$, and for each $n\geqslant 1$, we choose a point $x_n$ which is the generic point of an irreducible component $X_n'\subseteq (X_n\cap\overline{\{x_{n-1}\}})$ such that $X_n'\cap X_{n+1}\neq\emptyset$. 
		As a consequence, applying the dimension function to the sequence $\{x_n\}_{n=0,\ldots,m}$, we obtain a strictly descending sequence of integers $\{\delta(x_n)\}_{n=0,\ldots,m}$. However, since $\delta$ is bounded, we get a limit on the length $m$ of the descending chain $\{X_n\}$, implying that $|X|$ is of finite Krull dimension. %to be arbitrarily large contradicting the fact that $\delta$ is bounded. A Noetherian topological space has finite Krull dimension.  Since the fibers are Noetherian then it is Noetherian. Given a sequence $X_n$ fibrewise it is constant and there are finite number of finite (finite rank). Thus we are done.
		\par Hence, it suffices to show that $\delta$ is a bounded dimension function. Consider a specialisation $x\rightsquigarrow y$ in $X$. If $f(x)=f(y)$, then replacing $X$ by its fibre over $f(x)$, we may assume that $X$ is a finite type $\kappa(f(x))$-scheme; in which case, thanks to \cite[\href{https://stacks.math.columbia.edu/tag/02JW}{Tag 02JW}]{stacks-project}, $\delta(x)\geqslant\delta(y)$ and the specialisation is immediate if and only if $\delta(x)=\delta(y)+1$. Henceforth, we assume that $f(x)\neq f(y)$. Localising at $f(y)$, the function \bud\delta'\colon |X_{R_{f(y)}}|\to\bb{Z}\text{, given by, }\delta'(x)=\trdeg_{\kappa(f(x))}(\kappa(x))-\codim_{\spec R_{f(y)}}(\overline{\{f(x)\}})\eud equals $\delta|_{X_{R_{f(y)}}}\colon |X_{R_{f(y)}}|\to\bb{Z}$, up to a constant. Thus, localising at the prime ideal corresponding to $f(y)$, without loss of generality, we might assume that $R$ is a valuation ring with closed point $f(y)$. Further, dividing by the prime ideal corresponding to $f(x)$, we may also assume that $f(x)$ is the generic point. Therefore, the closed subscheme $Z\colonequals\overline{\{x\}}\subseteq X$ (with reduced structure) is dominant, producing the equality $\dim Z_{f(y)}=\dim Z_{f(x)}$ thanks to Lemma~\ref{lem:valuation_ring_fibres_same_dimension}. Moreover, since $Z$ is a dominant, integral $R$-scheme, it is automatically $R$-flat (it follows from the fact that flatness can be checked locally and from \ccite{bourbaki}{Chapter~I, \textsection2.4, Proposition~3(ii)}, which implies that an injection $R\hookrightarrow A$ into an integral domain is flat); additionally, since $Z$ is of $R$-finite type, by \ccite{raynaud_gruson}{premi\`ere partie, corollaire~3.4.7}, it is of $R$-finite presentation. On the other hand, since $x$ is the generic point of $Z$, it is also the generic point of $Z_{\kappa(x)}$; consequently, applying, for example Noether normalisation \cite[\href{https://stacks.math.columbia.edu/tag/00P0}{Tag 00P0}]{stacks-project}, we deduce that \bud\textstyle\trdeg_{\kappa(f(x))}(\kappa(x))=\dim Z_{f(x)}=\dim Z_{f(y)}\geqslant \trdeg_{\kappa(f(y))}(\kappa(y)).\eud Finally, since $\overline{\{f(x)\}}\supsetneq\overline{\{f(y)\}}$, the inequality $\delta(x)>\delta(y)$ follows. 
		\par Lastly, we show that the specialisation $x\rightsquigarrow y$ is immediate if and only if $\delta(x)=\delta(y)+1$. In similar vein as before, $\trdeg_{\kappa(f(x))}(\kappa(x))\geqslant\trdeg_{\kappa(f(y))}(\kappa(y))$, with equality in the case that $y$ is the generic point of $Z_{f(y)}$. As a consequence, $\delta(x)=\delta(y)+1$ is equivalent to the case when $\dim(\overline{\{f(x)\}})=\dim(\overline{\{f(y)\}})+1$ and $\trdeg_{\kappa(f(x))}(\kappa(x))=\trdeg_{\kappa(f(y))}(\kappa(y))$, which in turn is equivalent to the case when $f(y)$ is an immediate specialisation of $f(x)$ and $y$ is the generic point of $Z_{f(y)}$, in other words, $y$ is an immediate specialisation of $x$. Indeed, if $y$ is an immediate specialisation of $x$, then, $f(y)$ is an immediate specialisation of $f(x)$ (see \cite[\href{https://stacks.math.columbia.edu/tag/0D4H}{Tag 0D4H}]{stacks-project}).
		\par To verify that $\delta$ is bounded, we choose the generic point $x\in X$ of an irreducible component. Let $y\in\overline{\{x\}}$ be a closed point. Then, $0\leqslant \delta(x)-\delta(y)\leqslant\trdeg_{\kappa(f(x))}(\kappa(x))+\dim(\overline{\{f(x)\}})=\dim Z_{f(x)}+\dim(\overline{\{f(x)\}})\leqslant \dim X_{f(x)}+\dim R$, and hence we are done.    
	\es
	\br\label{rem:valuation_rings_catenary}
	Moreover, in Proposition~\ref{lem:valuation_rings_catenary}, if $R$ is semilocal, then $|Y|$ is even noetherian. Indeed, since $\spec(R)$ is then a finite set (see, for example, \ccite{amar:unramified_case_gersten_conjecture}{Remark~2.7}), this follows by noting that each of the $R$-fibres is noetherian.
	\er
	%While working with smooth algebras $A$ over a valuation ring $V$, to use local arguments, it is often useful to localise at a generic point $x$ in the $V$-special fibre of $\spec A$ (see \Cref{lem:Brauer_group_injects} and \Cref{prop:reduction_to_relative_curve_assuming_Borel}). The following lemma shows that $A_x$ is, in fact, a valuation ring.
%%%%%%%%%%%%%%%%%%%%%%%%%%%%%%%%%%%%%%%%%%%%%%%%%%%%%%%%%%%%%%%%%%%%%%%%%%%%%%%%%%%%%%%%%%%%%%%%%%%%%%%%%%%%%%%%%%%%%%%%%%%%%%%%%%%%%%%%%%%%%%%%%%%%%%%%%%%%%%%%%%%%%%%%%%%%%%%%%%%%%%%%%%%%%%%%%%%%%%%%%%%%%%%%%%%%%%%%%%%%%%%%%%%%%%%%%%%%%%%%%%%%%%%%%%%%%%%%%%%%%%%%%%%%%%%%%%%%%%%%%%%%%%%%%%%%%%%%%%%%%%%%%%%%%%%%%%%%%%%%%%%%%%%%%%%%%%%%%%%%%%%%%
%\section{Gabber's presentation lemma over mixed-characteristic valuation rings.}\label{section:presentation_lemma}
	 While working with non-noetherian rings, it is crucial to distinguish between finite type objects and finitely presented ones. This distinction becomes especially important when using noetherian approximation techniques, where finitely presented are required, and not merely finite type. Coherent rings (Definition~\ref{defn:coherent_rings}) form a significant class of rings where the gap between finite type and finite presentation is reduced.
	 \par We first need to define the notion of coherent modules, which we do below. %Thankfully, rings in $\mathcal{S}_{\mathrm{val}}$ are coherent, which permits us to adapt some of the techniques from \cite{gillet_levine}.  
	 \begin{montobo}[Coherence]\label{montobo:locally_coherent}
	 	Given a locally ringed space $X$, an $\ca{O}_X$-module $\scr{F}$ is called \textit{coherent} if it is of finite type and for every open $U\subseteq X$ and every finite collection $s_i\in\scr{F}(U)$, $i=1,\ldots,n$, the kernel of the associated morphism $\bigoplus_{i=1,\ldots,n}\ca{O}_U\to\scr{F}$ is of finite type (\cite[\href{https://stacks.math.columbia.edu/tag/01BV}{Tag 01BV}]{stacks-project}). 
	 \end{montobo}
	We note some of the properties of coherence below.
	 	\bn
	 		\item A coherent $\ca{O}_X$-module is finitely presented, and therefore, quasi-coherent (\cite[\href{https://stacks.math.columbia.edu/tag/01BW}{Tag 01BW}]{stacks-project}). 
	 		\item A finite type $\ca{O}_X$-submodule of a coherent $\ca{O}_X$-module $\scr{F}$ is coherent, and the same holds for a finite type $\ca{O}_X$-module that is a quotient of $\scr{F}$. 
	 		%\item A finite type $O_x$ module is coherent. A module is coherent iff it is coherent
	 		\item\label{part:weak_Serre_subcat} Furthermore, the category of coherent $\ca{O}_X$-modules form a weak Serre subcategory of the quasi-coherent $\ca{O}_X$-modules (see \stacks{01by}). In particular, this category is abelian.
	 	\en 
	 	\br
	 		It is challenging to work with coherent modules over arbitrary schemes. However, over a coherent scheme (Definition~\ref{defn:coherent_rings}), they are identical to the smallest abelian subcategory generated by the locally free sheaves of finite rank (see Remark~\ref{rem:coherence_finite_type_modules}).
	 	\er
	 	\bdf\label{defn:coherent_rings} A scheme $X$ is called \textit{locally coherent} if $\ca{O}_X$ is a coherent module over itself (\ccite{gabber-ramero-foundations}{Definition~8.1.54}) and $X$ is called \textit{coherent} (resp., \textit{stably coherent}) if it is locally coherent, quasi-compact and quasi-separated (resp., any $X$-scheme of finite presentation is coherent). A ring $A$ is called \textit{coherent} (resp., \textit{stably coherent}) if $\spec(A)$ is a coherent scheme (resp., a stably coherent scheme). 
	 	\edf
	 	\br\label{rem:coherence_finite_type_modules}
	 		On a locally coherent scheme $X$, a quasi-coherent $\ca{O}_X$-module of finite type is coherent. Indeed, this follows from the definitions.
	 	\er
	 	A determining property of a coherent ring is the following. A ring $A$ is coherent if and only if any finitely generated ideal $I\subseteq A$ is finitely presented (\cite[\href{https://stacks.math.columbia.edu/tag/05CV}{Tag 05CV}]{stacks-project}). For example, any noetherian ring is coherent, and as a consequence, stably coherent. An important class of non-noetherian rings that are stably coherent are Prüfer domains (see Lemma~\ref{lem:algebraic_space_coherent}), which we introduce below.%
	 	\bdf[\ccite{gilmer_multiplicative_ideal_theory}{Chapter~IV, Section~22}]\label{defn:prufer_domain}
	 		A commutative ring is said to be a {\it Prüfer domain} if it is an integral domain whose localisation at every prime ideal is a valuation ring. A commutative ring is said to be a {\it Prüfer ring} if it is a finite product of Prüfer domains.
	 		%\begin{enumerate}
	 		%    \item A {\it Prüfer domain} is an integral domain whose localisation at every prime ideal is a valuation ring.
	 		%    \item A {\it Prüfer ring} is a finite product of Prüfer domains.
	 		%\end{enumerate}
	 	\edf
	 	
	 	\br\label{rem:Prufer_ring} 
	 		\bn[(i)] 
	 			\item\label{rem:prufer_ring_pt:1} A Prüfer ring, being a product of integral domain, is reduced. In particular, a connected Prüfer ring is integral, equivalently, a Prüfer domain.
	 			\item\label{rem:prufer_ring_pt:2} A Prüfer ring of Krull dimension $0$ is a finite product of fields.
	 			\item Valuation rings themselves are Prüfer domains, and, in fact, any local Prüfer ring must be a valuation ring. As an example of a non-local Prüfer domain, we have the ring of algebraic integers $\overline{\bb{Z}}$.
	 			\item The class of Dedekind rings is equal to the class of noetherian Prüfer rings. However, contrary to Dedekind rings, Prüfer rings can be of arbitrary Krull dimension. For instance, the subring $\{P(X) \in \bb{Q}[X] \text{ } \vert \text{ } P(0) \in \bb{Z}\}$ of $\bb{Q}[X]$ is a Prüfer ring of Krull dimension two (\ccite{cahen_chabert_integer_valued_polynomials}{Theorem 17}), and the ring of entire holomorphic functions on the complex plane is a Prüfer ring of infinite Krull dimension (\cite{loped_entire_functions}).
	 			\item\label{rem:prufer_ring_pt:3} An integral domain is a Prüfer domain if and only if each of its nonzero finitely generated ideals is invertible (\ccite{gilmer_multiplicative_ideal_theory}{Theorem~$22.1$}). %\st{As noted in op.~cit.~Theorem 22.1, a ring $R$ is a Prüfer domain if and only if each nonzero finitely generated ideal is invertible.} 
	 			By \ccite{bourbaki}{Chapter I, 2.4, Prop. 3(ii)}, this implies that a module over a Prüfer ring is flat if and only if it is torsionfree.
	 		\en
	 	\er
%	 	\bl[\stacks{07ub}]
%	 	Let $X$ be a locally coherent algebraic stack. Given a quasi-coherent $\ca{O}_X$-module $\scr{F}$, it is coherent in either of the following cases.
%		 	\bn[(a)]
%		 	\item There exists a surjective étale morphism $\varphi\colon U\to X$ from an algebraic space so that $\varphi^{\ast}\scr{F}$ is a coherent module on $U$.
%		 	\item For any surjective étale morphism $\varphi\colon U\to X$ from an algebraic space, the pullback $\varphi^{\ast}\scr{F}$ is a coherent module on $U$.
%		 	\item $\scr{F}$ is a finite type $\ca{O}_X$-module.
%	 	%\item $\scr{F}$ is a finitely presented $\ca{O}_X$-module.
%	 	\en
%	 	\el
	 	\bl\label{lem:algebraic_space_coherent}% 
	 		A scheme locally of finite presentation over a Pr\"ufer ring is locally coherent.
	 	\el
	 	\bs Indeed, since the property of being locally coherent is étale-local (see \stacks{05vr}), it suffices to check that any ring $A$ that is a finitely presented algebra over a Prüfer ring $R$ is coherent. Let $f\colon A'\colonequals R[x_1,\ldots,x_n]\twoheadrightarrow A$ be a presentation of $A$ such that $\ker(f) \subset A'$ is a finitely generated ideal. Since $\ker(f) \subset A'$ is finitely generated, it is enough to show that the ring $A'$ is coherent. Letting $I\subset A'$ be a finitely generated ideal, we shall show that $I$ is a finitely presented $A'$-module. Putting $X=\spec A', S=\spec R$ and $\scr{M}=\widetilde{I}$ in \ccite{raynaud_gruson}{premi\`ere partie, théorème~3.4.6} (by \ccite{bourbaki}{Chapter~I, \textsection2.4, Proposition~3(ii)}, \cite[\href{https://stacks.math.columbia.edu/tag/090Q}{Tag 090Q}]{stacks-project} and the fact that flatness is a local property \cite[\href{https://stacks.math.columbia.edu/tag/0250}{Tag 0250}]{stacks-project}, the $R$-torsion-free module $I$ is flat), we obtain that $I$ is a finitely presented $A'$-module, showing that $A'$ is coherent.\es
	 	%The schemes considered in this chapter are flat and finite type over Pr\"ufer domains. Let $X$ be such a scheme; more precisely, let $X$ be a flat, finite type scheme over a Pr\"ufer domain $R$. Thanks to \ccite{raynaud_gruson}{Premi\`ere partie, Corollaire~3.4.7}, $X$ is $R$-finitely presented, and hence, by the discussion above, it is a coherent scheme. As a result, the schemes considered in this chapter are coherent.
	 We are now prepared to show that higher direct images along projective morphisms preserve coherence of sheaves. Our argument follows \ccite{egaIII_1}{théorème~2.2.1 and corollaire~2.2.4}, and generalises \ccite{phd-thesis}{Lemma~6.2}.
	 \bl\label{lem:ega_III_1_2.2.2}
		 Let $Y$ be a stably coherent scheme and and let $f\colon X\to Y$ be a projective morphism of finite presentation with a closed immersion $\iota\colon X\hookrightarrow\bb{P}^{m-1}_Y$, for some $m\geqslant 1$. Let $\ca{O}_X(1)\colonequals\iota^{\ast}(\ca{O}_{\bb{P}^{m-1}_Y}(1))$ and for any quasi-coherent $\ca{O}_X$-module $\scr{G}$, let $\scr{G}(n)\colonequals\scr{G}\otimes_{\ca{O}_X}\ca{O}_X(n)$. Then, given a surjection $\varphi\colon\scr{F}\twoheadrightarrow\scr{G}$ of coherent $\ca{O}_X$-modules,
		 \bn[(i)]
			 \item\label{part:1} we have $R^{q}f_{\ast}\scr{F}=0$, for each $q>m$, 
			 \item\label{part:2} there exists an integer $N$ such that for all $n\geqslant N$, we have \bud\text{$R^{q}f_{\ast}(\scr{F}(n))=0$ for any $q\geqslant 1$,}\eud
			 \item\label{part:4} the $\ca{O}_Y$-module $R^qf_{\ast}(\scr{F})$ is coherent for any $q$, and
			 \item\label{part:3} there exists an integer $N$ such that for all $n\geqslant N$, we have \bud f_{\ast}(\varphi)\colon f_{\ast}(\scr{F}(n))\twoheadrightarrow f_{\ast}(\scr{G}(n)).\eud
		 \en
	 \el
	 \bs
	 Since the statements are Zariski local on the base, we may assume that $Y=\spec(R)$, where $R$ is a stably coherent ring. Since both $X$ and $\bb{P}^{m-1}_R$ are of $R$-finite presentation, they are coherent schemes.
	 %We note that $X, Y$ and $\bb{P}^{m-1}_Y$ are actually of finite presentation over $S$, since they are $S$-flat of finite type (see \ccite{raynaud_gruson}{Première partie, théorème~3.4.6}). In particular, they are coherent algebraic spaces.
	 %Cech complex for algebraic spaces \stacks{0721}, proper base change of algebraic spaces \stacks{0A1N}.
	 %Given a quasi-coherent $\ca{O}_X$-module $\scr{K}$, by \ccite{egaII}{Corollaire~3.4.5 and Proposition~3.5.2}, we have $\iota_{\ast}(\scr{K}(n))\cong(\iota_{\ast}\scr{K})(n)$, for all $n\geqslant 0$. 
	 \vspace{0.25cm}
	 \par\eqref{part:1}$\colon$ Since $X$ is a closed subscheme of $\bb{P}^{m-1}_R$, it can be covered by $m$ affines, say $\{U_i\}$. Consequently, thanks to \cite[\href{https://stacks.math.columbia.edu/tag/01XD}{Tag 01XD}]{stacks-project} or \ccite{egaIII_1}{proposition~1.4.1}, since $X$ is separated, the $q$-th Čech cohomology group $\check{H}^q(\{U_i\},\scr{F})$ of $\scr{F}$ with respect to $\{U_i\}$ identifies itself with $H^q(X,\scr{F})$, for each $q$. The claim follows because $\check{H}^q(\{U_i\},\scr{F})=0$, for all $q>m$. \vspace{0.25cm} %(\cite[\href{https://stacks.math.columbia.edu/tag/02UZ}{Tag 02UZ}]{stacks-project}).
	 \par\eqref{part:2}$\colon$ We follow the proof of \ccite{egaIII_1}{proposition~2.2.2}. Given that $\bb{P}^{m-1}_R$ is coherent, this means that $\iota_{\ast}(\ca{O}_X)$ is automatically a coherent $\ca{O}_{\bb{P}^{m-1}_R}$-module (\cite[\href{https://stacks.math.columbia.edu/tag/01BZ}{Tag 01BZ}]{stacks-project}). Similarly, $\iota_{\ast}(\scr{F})$ is a coherent $\ca{O}_{\bb{P}^{m-1}_R}$-module. Since higher direct images under a closed immersion vanish (\cite[\href{https://stacks.math.columbia.edu/tag/01QY}{Tag 01QY}]{stacks-project}), it is enough to show that exists an integer $N$ such that for any $n\geqslant N$, we get  $H^{q}(\bb{P}^{m-1}_R,\iota_{\ast}(\scr{F}(n)))=0$, for all $q\geqslant 1$. Therefore, it suffices to assume that $X=\bb{P}^{m-1}_R$. 
	 \par Thanks to \ccite{egaII}{corollaire~2.7.10}, there exists a surjection $j\colon\scr{L}\twoheadrightarrow\scr{F}$, where $\scr{L}$ is a finite direct sum of modules of the form $\ca{O}_X(r)^{\oplus s}$ for some $r\in\bb{Z}$ and $s\geqslant 0$. Letting $\scr{K}\colonequals\ker j$, we get a short exact sequence $0\to\scr{K}\to\scr{L}\to\scr{F}\to 0$ of coherent sheaves (\cite[\href{https://stacks.math.columbia.edu/tag/01BY}{Tag 01BY}]{stacks-project}). Since $\ca{O}_X(n)$ is a locally free $\ca{O}_X$-module for any $n$, we get a short exact sequence \bd\label{short_exact:twisitng}\text{$0\to\scr{K}(n)\to\scr{L}(n)\to\scr{F}(n)\to 0$},\ed of coherent modules (\cite[\href{https://stacks.math.columbia.edu/tag/01CE}{Tag 01CE}]{stacks-project}), for each $n$. We shall show \eqref{part:2} by the method of descending induction. For $q>m$, the result follows from \eqref{part:1}. Suppose that for $d\geqslant2$ and for any coherent $\ca{O}_X$-module $\scr{M}$, there exists an integer $N$ such that $H^{q}(X,\scr{M}(n))=0$, for all $q\geqslant d$ and for any $n\geqslant N$. We shall show that there exists an integer $N$ such that $H^{q}(X,\scr{F}(n))=0$, for all $q\geqslant d-1$ and for any $n\geqslant N$. Thanks to \ccite{egaIII_1}{corollaire~2.1.13}, we have $H^{q}(\ca{O}_X(t))=0$, for all $q\geqslant 1$ and for any $t\geqslant 0$; consequently, \bud\text{$H^{q}(X,\scr{L}(n))=0$, for all $q\geqslant 1$ and $n\gg 0$}.\eud We choose $N$ such that for any $n\geqslant N$, we get $H^{q}(X,\scr{L}(n))=0$, for every $q\geqslant 1$, and $H^{q}(X,\scr{K}(n))=0$, for every $q\geqslant d$. With this choice, writing the associated long exact sequence of cohomology of \eqref{short_exact:twisitng}, we get isomorphisms $H^{q}(X,\scr{F}(n))\cong H^{q+1}(X,\scr{K}(n))$, for all $q\geqslant 1$ and for any $n\geqslant N$. This implies that $H^{q}(X,\scr{F}(n))=0$, for every $q\geqslant d-1$ and for any $n\geqslant N$, and the induction step is complete. Thus, we are done. \vspace{0.25cm}
	 \par\eqref{part:4}$\colon$ Our proof by descending induction shall be similar to \eqref{part:2}. We consider the long exact sequence of cohomology \bud \cdots\to H^{q-1}(X,\scr{K})\to H^{q-1}(X,\scr{L})\to H^{q-1}(X,\scr{F})\to H^q(X,\scr{K})\to H^q(X,\scr{L})\to\cdots\eud associated to \eqref{short_exact:twisitng}. Thanks to \eqref{part:weak_Serre_subcat}, we reduce to establish the claim for $\scr{F}\cong\ca{O}_X(r)$ for each $r$. In this case, it follows from the nature of cohomology of projective spaces \ccite{egaIII_1}{proposition~2.1.12}. \vspace{0.25cm}
	 \par\eqref{part:3}$\colon$ Letting $\scr{K}\colonequals\ker\varphi$, we get a short exact sequence $0\to\scr{K}\to\scr{F}\to\scr{G}\to 0$ of coherent sheaves (\cite[\href{https://stacks.math.columbia.edu/tag/01BY}{Tag 01BY}]{stacks-project}). In a similar vein as above, since $\ca{O}_X(n)$ is a locally free $\ca{O}_X$-module for any $n$, we get a short exact sequence \bd\label{short_exact:twisitngII}\text{$0\to\scr{K}(n)\to\scr{F}(n)\to\scr{G}(n)\to 0$}\ed of coherent modules (\cite[\href{https://stacks.math.columbia.edu/tag/01CE}{Tag 01CE}]{stacks-project}), for each $n$. By \eqref{part:2}, there exists an integer $N$ such that for any $n\geqslant N$, we have $H^1(X,\scr{K}(n))=0$. Writing the long exact sequence of cohomology associated to \eqref{short_exact:twisitngII}, we get the requisite surjection, and we are done.
	 \es
	 The remainder of the section is dedicated to the proof of our main result, Theorem~\ref{prop:moving_lemma_general}, which is a version of Gabber's presentation lemma over relatively general base rings, including arbitrary noetherian domains and Prüfer rings. Both the statement, as well as the proof, are inspired by, and extend, \ccite{ces_grothendieck-serre}{Variant 3.7}, which treats the case when the base is a Dedekind ring. At the same time, our result generalises \ccite{phd-thesis}{Proposition~6.4}, which establishes the case when the base is a valuation ring of finite rank. 
	 \par Our approach proceeds by base changing to the special fibres and then bootstrapping from a presentation obtained there. Specifically, we apply \ccite{ces_grothendieck-serre}{Proposition~3.6}, which furnishes such a presentation when the base is a field.
	 \par Consequently, it is essential--especially when $n>1$--that the points $x_1,\ldots,x_n$ specialise to the special fibres. The following result allows us reduce to this situation in the proof of Theorem~\ref{prop:moving_lemma_general}. 
	 \bl\label{lem:specialises_to_closes_fibre}
	 Let $R$ be a semilocal Pr\"ufer domain of finite Krull dimension, let $X$ be a flat, projective $R$-scheme that is $R$-fibrewise of pure dimension $d$, let $\ca{O}_X(1)$ be an $R$-relatively very ample line bundle on $X$, let $W\subseteq X^{\mathrm{sm}}$ be an open, and let $Y\subset X$ be a closed subscheme such that $Y\setminus W$ is $R$-fibrewise of codimension $\geqslant 2$. Then, given any points $x_1,\ldots,x_n\in W$, there exist
	 \bun
	 \item[$\circ$] a semilocal Pr\"ufer domain $\tilde{R}$ of finite Krull dimension with an open subset $\spec(R)\subseteq\spec(\tilde{R})$, 
	 \item[$\circ$] a flat, projective $\tilde{R}$-scheme $\tilde{X}$ which is $\tilde{R}$-fibrewise of pure dimension $d$ extending the $R$-scheme $X$,
	 \item[$\circ$] an $\tilde{R}$-relatively very ample line bundle $\ca{O}_{\tilde{X}}(1)$ on $\tilde{X}$ whose restriction to $X$ is $\ca{O}_X(1)$,
	 \item[$\circ$] an open $\tilde{W}\subseteq \tilde{X}^{\mathrm{sm}}$ whose intersection with $X$ is $W$, and
	 \item[$\circ$] a closed subscheme $\tilde{Y}\subset\tilde{X}$ whose restriction to $X$ is $Y$ such that $\tilde{Y}\setminus\tilde{W}$ is $\tilde{R}$-fibrewise of codimension $\geqslant 2$; 
	 \eun
	 so that each $x_1,\ldots,x_n$ specialises to an $\tilde{R}$-special fibre of $\tilde{W}$.
	 \el
	 \bs
	 	If each $x_i$ already specialises to a point in an $R$-special fibre of $W$--which includes the case when $R$ is a product of fields--then there is nothing to show. Otherwise, let $y_1,\ldots,y_m\in X$ be the points that fail to satisfy this condition, and $\ca{P}\subset\spec(R)$ denote their images.
	 	
	 	\par %In general, $x_1,\ldots,x_n$ might not have a specialisation in a special fibre $W_C$ of $W$, say $y_1,\ldots,y_m$. Let $\mathcal{P}\subseteq\spec(R)$ be the images of $y_1,\ldots,y_m$. 
	 	%The rest of the proof is similar to the same of \ccite{amar:unramified_case_gersten_conjecture}{Presentation Lemma~3.2}. 
	 	Proceeding as in the proof of \ccite{amar:unramified_case_gersten_conjecture}{Presentation Lemma~3.2}, %we shall tailor $X$ in such a way that each $y_i$ specialises to $W_C$, effectively replacing the original $X$ with this customised version. In this regard, 
	 	thanks to op.~cit.~Lemma~2.11(iii) and a limit argument, without loss of generality, we may assume that the residue field $\kappa(\fr{p})$ of $R$ at each $\fr{p}\in\cal{P}$ is finitely generated over its prime subfield. By, for example, \ccite{phd-thesis}{Lemma~6.1}, each field $\kappa(\fr{p})$ is a fraction field of a regular domain $A_{\fr{p}}$ that is smooth over $\bb{F}_p$ or $\bb{Z}$. Moreover, each $A_{\fr{p}}$ is of positive Krull dimension, since otherwise $K$ is a finite field, in which case, it contradicts our assumption that $R$ is not a field. By localising $A_{\fr{p}}$, we may assume that
	 	\bn[(1)]
	 	\item the scheme $X_{\kappa(\fr{p})}$ spreads out to a projective, flat $A_{\fr{p}}$-scheme $X_{\fr{p}}$ that is fibrewise of pure dimension $d$ by \ccite{egaIV_3}{théorème~12.2.1 (ii) and (v)},
	 	\item  the relative $\kappa(\fr{p})$-very ample line bundle $\ca{O}_{X_{\kappa(\fr{p})}}(1)$ spreads out to a relative $A_{\fr{p}}$-very ample line bundle,
	 	\item there is an open $W_{\fr{p}}\subset X_{\fr{p}}^{\text{sm}}$ which intersects the $\kappa(\fr{p})$-fibre at $W_{\kappa(\fr{p})}$ by \cite[\href{https://stacks.math.columbia.edu/tag/01V9}{Tag 01V9}]{stacks-project},
	 	\item each point $y_i$ that lies in $W_{\kappa(\fr{p})}$ spreads out to an $A_{\fr{p}}$-finite closed subscheme in $W_{\fr{p}}$,
	 	\item and the closed subscheme $Y_{\kappa(\fr{p})}$ spreads out to a closed subscheme $Y_{\fr{p}}\subset X_{\fr{p}}$ such that $Y_{\fr{p}}\setminus W_{\fr{p}}$ is $A_{\fr{p}}$-fibrewise of codimension $\geqslant 2$ in $X_{\fr{p}}$ (see \ccite{egaIV_3}{corollaire~12.2.2 (i)}).
	 	\en 
	 	Given that $A_{\fr{p}}$ is of positive Krull dimension, it has infinitely many primes of height $1$, allowing us to choose such a prime $\fr{r}\subset A$ so that the localisation $A'_{\fr{p}}$, which is necessarily a discrete valuation ring, of $A_{\fr{p}}$ at $\fr{r}$ is different from each of the localisations of $R/\fr{p}$. Choosing such a prime $\fr{r}\subset A$, we substitute $A_{\fr{p}}$ with $A'_{\fr{p}}$ and consider the affine scheme $\spec(\tilde{R})$ obtained by gluing $\spec(A_{\fr{p}})$ to $\spec(R)$ at $\spec(\kappa(\fr{p}))$. Thanks to \ccite{amar:unramified_case_gersten_conjecture}{Lemma~2.6}, the resulting ring $\tilde{R}$ is a semilocal Pr\"ufer domain of finite Krull dimension. Similarly,
	 	\bn [(1)]
	 	\item we glue $X$ and $X_{\fr{p}}$ along $X_{\kappa(\fr{p})}$ to obtain a projective $\tilde{R}$-scheme $X_{\tilde{R}}$ fibrewise of pure dimension $d$ with a relative very ample line bundle,  
	 	\item we glue $W$ and $W_{\fr{p}}$ along $W_{\kappa(\fr{p})}$ to obtain an open $W_{\tilde{R}}\subset X_{\tilde{R}}^{\text{sm}}$, 
	 	\item we glue $Y$ and $Y_{\fr{p}}$ along $Y_{\kappa(\fr{p})}$ to obtain a closed subscheme $Y_{\tilde{R}}\subset X_{\tilde{R}}$ such that the special fibres of $Y_{\tilde{R}}\setminus W_{\tilde{R}}$ are of codimension $\geqslant 2$.
	 	\en
	 	By systematically advancing through the primes $\fr{p}\in\cal{P}$, ordered by their height, we can gradually build $\tilde{R}$ and the corresponding objects as described above. As a result of this construction, even $y_1,\ldots,y_m$ specialise to points in the special fibre of $W_{\tilde{R}}$, finishing the proof. %Consequently, we can apply the previous case to complete our proof.
	 \es
%	 \bdf\label{defn:noetherian_like}
%	 	\rl{An integral domain $R$ is said to be \textit{noetherian-like} if it is stably coherent and the topological space associated with any finite type scheme over a semilocalisation of $R$ at finitely many primes is noetherian.}
%	 \edf
%	 \br
%	 	\bun 
%	 		\item[$\circ$] A noetherian domain is noetherian-like. 
%	 		\item[$\circ$] The same holds for any Prüfer domain of finite Krull dimension. Indeed, any Prüfer ring is stably coherent and the topological space associated with any finite type scheme over a semilocal Prüfer ring of finite Krull dimension is noetherian (see \Cref{rem:valuation_rings_catenary}).
%	 	\eun
%	 \er
	 Let us prove our main theorem below.
	 \bt\label{prop:moving_lemma_general}
	 Let $R$ be ring and let $X$ be a projective $R$-scheme of finite presentation. Let $\ca{O}_X(1)$ be an $R$-relatively very ample line bundle on $X$, and let $W\subseteq X^{\mathrm{sm}}$ be a quasi-compact open subset that is $R$-fibrewise of dimension $d>0$ and contains points $x_1,\ldots,x_n$, for some integer $n\geqslant 1$. Let $Y\subset X$ be a finitely presented closed subscheme such that $Y\setminus W$ is $R$-fibrewise of codimension $\geqslant 2$. We assume further that:
	 \bn[(a)] 
	 \item\label{case:a} either $R$ is a Prüfer ring, or
	 \item\label{case:b} $n=1$ and $R$ is a stably coherent ring in the sense of Definition~\ref{defn:coherent_rings} (for example, a noetherian ring).
	 \en 
	 Then, letting $w_1\colonequals 1$, after replacing $\ca{O}_X(1)$ by a sufficiently large power, there exist integers $w_2,\ldots,w_d\geqslant 1$ and nonzero sections $s_k\in\Gamma(X, \ca{O}(w_k))$ for each $k=1,\ldots,d$, as well as, affine opens \bud S\subseteq\bb{A}^{d-1}_{R}\text{\hspace{0.75cm}and\hspace{0.85cm}}x_1,\ldots,x_n\in U\subseteq W\cap\pi^{-1}(S)\eud such that the morphism $\pi\colon U\to S$, determined by the sections $s_i$, is smooth of relative dimension $1$ and $Y\cap U=Y\cap\pi^{-1}(S)$ is $\pi$-finite. Moreover, if $W$ is $R$-fibrewise of pure constant dimension, then $\pi$ is automatically of pure relative dimension $1$. 
%	 \bn[(i)]
%	 \item
%	 \item
%	 \item
%	 \item\label{condition:iv} for every $i$, we have $Y\cap H_1\cap\overline{\pi}^{-1}(\pi(x_i))=\emptyset$,
%	 \item\label{condition:v} for every $i$, the morphism $\pi$ is smooth at each $Y\cap\overline{\pi}^{-1}(\pi(x_i))$, and
%	 \item\label{condition:vi} there are affine opens \bud S\subseteq\bb{A}^{d-1}_{R}\text{\hspace{0.75cm}and\hspace{0.85cm}}x_1,\ldots,x_n\in U\subseteq W\cap\pi^{-1}(S)\eud such that $\pi\colon U\to S$ is smooth of relative dimension $1$ and $Y\cap U=Y\cap\pi^{-1}(S)$ is $\pi$-finite.
	 %\en
	 % and affine opens $S\subset\bb{A}^{d-1}_R$ and $x_1,\ldots,x_n\in U\subset W$ such that $U$ intersects $Y\cap W$ in a nonempty set and that the morphism $\pi\colon U\to S$ determined by $h_1,\ldots,h_{d-1}$ is smooth of relative dimension $1$ and that its restriction to $Y\cap U$ is finite.
	 \et
	 \bs
	 %The proof of \eqref{point:i}-\eqref{point:vi} follow from the corresponding point in \Cref{prop:moving_lemma_field}. Indeed, we use the fact that a closed subset $Z\subset\overline{X}$ is empty if $Z\cap C$ is empty, for each $R$-special fibre $C$ of $\overline{X}$. 
	 If $R$ is a product of fields, then the claim follows from \ccite{ces_grothendieck-serre}{Proposition~3.6}. Therefore, we may assume that $R$ is not a product of fields. 
	 \par Restricting to a connected component of $\spec(R)$, in either case, we may further assume that $\spec(R)$ is connected. In particular, in case \eqref{case:a}, this $R$ is a Prüfer domain (see Remark~\ref{rem:Prufer_ring}\eqref{rem:prufer_ring_pt:1}) of positive Krull dimension (see Remark~\ref{rem:Prufer_ring}\eqref{rem:prufer_ring_pt:2}). Moreover, in the same case, it suffices to consider when $R$ is of finite Krull dimension. Indeed, \ccite{amar:unramified_case_gersten_conjecture}{Lemma~2.5(b)} proves that $R$ is an increasing union of its Prüfer sub-domains $R_{\lambda}$ of finite Krull dimension. These canonical maps $R_{\lambda}\hookrightarrow R$ are, in fact, flat (see Remark~\ref{rem:Prufer_ring}\eqref{rem:prufer_ring_pt:3}), and since flat morphisms preserve fibrewise dimension, we may descend all data to some $R_{\lambda}$ (thanks to \stacks{0ey2}, the fibrewise dimension of $X$ is preserved, in addition, by \stacks{0h3v}, the fibrewise codimension of $Y\setminus W$ is preserved), and then base change back to $R$.
	 %\par \rl{Since $R$ is a filtered union of Prüfer subdomains of finite Krull dimension, these transition maps are affine. Thus, by a limit argument, we may reduce to the case of a Prüfer domain of finite Krull dimension.}
	 \par In both cases, since the claim is Zariski-local around $x_1,\ldots,x_n\in X$, we may, without loss of generality, localise $R$ at the images of $x_1,\ldots,x_n$ and then ultimately spread out to assume that $R$ is semilocal (which is local if $n=1$). Let $C$ be the reduced subscheme of closed points of $\spec(R)$. When $n>1$, by our assumption $R$ is a semilocal Prüfer domain of finite Krull dimension, in which case, Lemma~\ref{lem:specialises_to_closes_fibre} demonstrates that each of $x_1,\ldots,x_n$ specialises to a point in $W_C$. This is, however, automatically true in the case $n=1$. Therefore, in either case, each of the points $x_1,\ldots,x_n$ specialises even to a closed point $x'_i$ in $W_C$. We shall, without loss of generality, specialise each point $x_i$ to $x'_i$ and assume that $x_i=x'_i$, i.e., we assume that $x_i$ is a closed point in $W_C$.
	 \par Let $I\subset R$ be the ideal of vanishing of $C$. We write $I=\bigcup I_{\lambda}$, where the filtered union is taken over the set of finitely generated sub-ideals $I_{\lambda}\subset I$, and set $C_{\lambda}\colonequals\spec(R/I_{\lambda})$, for each $\lambda$. Letting $\iota_{\lambda}\colon X_{C_{\lambda}}\colonequals X\times_{\spec R}C_{\lambda}\hookrightarrow X$ be the inclusion, since $I_{\lambda}$ is finitely generated, the $\ca{O}_X$-module $\iota_{\lambda\ast}\ca{O}_{X_{C_\lambda}}$ is finitely presented, and hence, the morphism $\ca{O}_X\to\iota_{\lambda\ast}\ca{O}_{X_{C_\lambda}}$ is a surjection of coherent $\ca{O}_X$-modules (Remark~\ref{rem:coherence_finite_type_modules} and Lemma~\ref{lem:algebraic_space_coherent}). As a consequence,  by Lemma~\ref{lem:ega_III_1_2.2.2}\eqref{part:3}, there exists an integer $N$ such that \bd\label{surjection}\text{$\Gamma(X,\ca{O}_X(r))\twoheadrightarrow\Gamma(X,(\iota_{\lambda\ast}\ca{O}_{X_{C_{\lambda}}})(r))=\Gamma(X_{C_{\lambda}}, \ca{O}_{X_{C_{\lambda}}}(r))$ is a surjection, for all $r\geqslant N$}.\ed Replacing $\ca{O}_X(1)$ by a sufficiently large power, without loss of generality, we may assume that $N=1$ in \eqref{surjection}. Since, by our assumption, the points $x_1,\ldots,x_n$ lie over $C$, we use \ccite{ces_grothendieck-serre}{Proposition~3.6} to find sections $h_i\in\ca{O}_{X_C}(w_i)$, for each $i$, (the last aspect of loc.~cit.~ensures that these $h_i$ may be chosen to have constant degrees on $C$) that satisfy the claim in loc~cit. By a limit argument, there exist a $\lambda$ and sections $h_{i,\lambda}\in\ca{O}_{X_{C_{\lambda}}}(w_i)$ that lift $h_i$, for all $i$. Finally, \eqref{surjection} implies that there exist sections $s_i\in\ca{O}_X(w_i)$ that lift $h_{i,\lambda}$, for each $i$. %Upon raising $\ca{O}_X(1)$ to the power $w_1$, without loss of generality, we may assume that $w_1=1$.
	 \par Since $Y$ and $H_i$ are closed subschemes of the projective $R$-scheme $X$, for all $i$, they are $R$-projective (see \ccite{egaII}{définition~5.5.2}), in particular, $R$-proper. As a consequence, their images along the respective structure morphisms to $\spec(R)$ are closed, whence we get that \bud\text{the vanishing locus $H_1\colonequals V(s_1)$ does not contain any $x_i$}\eud and the vanishing loci $H_i\colonequals V(s_i)$ satisfy $$Y\cap H_1\cap\ldots\cap H_d=\emptyset$$ from their respective counterparts in \ccite{ces_grothendieck-serre}{Proposition~3.6}. Letting $\pi$ be the morphism determined by the sections $s_i$, we have a commutative diagram
	 \bud X\setminus H_1\arrow[d,"\pi"]\arrow[r,hook]&X\setminus H_1\cap\ldots\cap H_d\arrow[r,hook]\arrow[d,"\overline{\pi}"]&\overline{X}\colonequals\mathrm{Bl}_X(s_1,\ldots,s_d)\arrow[d,"\overline{\pi}"]\\\bb{A}_{R}^{d-1}\arrow[r,hook]&\bb{P}_{R}(w_1,\ldots,w_d)\arrow[r,-,double equal sign distance, double]&\bb{P}_{R}(w_1,\ldots,w_d).\eud
	  We note that in the displayed diagram above the weighted blowup need not commute with base change to $C$, however, the formation of the morphism $\overline{\pi}\colon X\setminus (H_1\cap\ldots\cap H_d)\to\bb{P}_R(w_1,\ldots,w_d)$ does and this suffices for our purposes. Thanks to the fibrewise criterion of flatness \cite[\href{https://stacks.math.columbia.edu/tag/039C}{Tag 039C}]{stacks-project}, the morphism $\pi$ of finite presentation is flat at $x_i$, for each $i$, whence, the fibrewise criterion of smoothness \cite[\href{https://stacks.math.columbia.edu/tag/01V8}{Tag 01V8}]{stacks-project} ensures also that it is smooth at $x_i$, for each $i$.
	  We now prove that for every $i$, we have \bd\label{diag:intersection_empty} Y\cap H_1\cap\overline{\pi}^{-1}(\pi(x_i))=\emptyset.\ed To do so, we use the fact that $\overline{\pi}$ is proper to argue that $\pi(x)=\overline{\pi}(x)\in\bb{P}_R(w_1,\ldots,x_d)$ is a closed point, which implies $\overline{\pi}^{-1}(\pi(x))\subset\overline{X}$ is a closed subset. However, by our choice of the sections $s_i$, and since images of proper morphisms are closed, \eqref{diag:intersection_empty} holds because it is true after base changing to $C$. In a similar vein as above, the proof that $\pi$ is smooth at each $Y\cap\overline{\pi}^{-1}(\pi(x))$ follows from the fibrewise criterion of flatness \cite[\href{https://stacks.math.columbia.edu/tag/039C}{Tag 039C}]{stacks-project} followed by the fibrewise criterion of smoothness \cite[\href{https://stacks.math.columbia.edu/tag/01V8}{Tag 01V8}]{stacks-project}. 
	  \par It remains to produce affine open subsets $U$ and $S$. First, we claim that the morphism $\overline{\pi}$ when restricted to $Y\cap\overline{\pi}^{-1}(\pi(x_i))$ has finite $R$-fibres (and hence, by \cite[\href{https://stacks.math.columbia.edu/tag/02NH}{Tag 02NH}]{stacks-project}, it is quasi-finite), for each $i$. Combining \eqref{diag:intersection_empty} and that $H_1$ is a hypersurface, this claim is a consequence of Krull's principal ideal theorem. The openness of the quasi-finite locus (\cite[\href{https://stacks.math.columbia.edu/tag/01TI}{Tag 01TI}]{stacks-project}) implies that there exists an open subset $U_1\subseteq Y$ containing $Y\cap\pi^{-1}(\pi(x_i))$, for all $i$, such that $\pi|_{U_1}$ is quasi-finite. Since $Y$ is proper, taking any open subset \bud\text{$\pi(x_1),\ldots,\pi(x_n)\in S_0\subseteq(\bb{A}^{d-1}_R\setminus\pi(Y\setminus U_1))$},\eud we observe that $\pi|_{Y\cap\pi^{-1}(S_0)}$ is quasi-finite, which implies that it is even finite (\stacks{02OG}). We choose an affine open $\pi(x_1),\ldots,\pi(x_n)\in S_0\subseteq(\bb{A}^{d-1}_R\setminus\pi(Y\setminus U_1))$. By the definition of a smooth morphism \cite[\href{https://stacks.math.columbia.edu/tag/01V5}{Tag 01V5}]{stacks-project}, there exists an affine open \bud\text{$U_0\subseteq\pi^{-1}(S_0)\cap W$}\eud containing $x_1,\ldots,x_n$ and the points of $Y\cap\overline{\pi}^{-1}(\pi(x_i))$, for all $i$, such that $\pi|_{U_0}\colon U_0\to S_0$ is smooth. A dimension count shows that $\pi|_{U_0}$ is of relative dimension $1$. Finally, it remains to find affine opens $\pi(x_1),\ldots,\pi(x_n)\in S\subseteq S_0$ and $U\subseteq U_0$ containing $x_1,\ldots,x_n$ and $Y\cap\overline{\pi}^{-1}(\pi(x_i))$, for all $i$, such that $Y\cap U=Y\cap\pi^{-1}(S)$. 
	  For this, we can choose any principal affine open 
	  \bud\text{$\pi(x_1),\ldots,\pi(x_n)\in S\subseteq S_0\setminus\pi(Y\setminus U_0)$ and set $U\colonequals U_0\cap\pi^{-1}(S)$}.\eud 
	  Finally, if $W$ is $R$-fibrewise of pure constant dimension, the same holds for $U$; consequently, the dimension count even shows us that $\pi|_U$ is of pure constant relative dimension. Hence, we are done.
	 %\vspace{1cm}
	 %\par Let $x\in Y\cap\overline{\pi}^{-1}(\pi(x_i))$ be a point. so that there are opens $S\subset\bb{A}^{d-1}$ and $x_1,\ldots,x_n\in U\subset W$ such that the morphism $\pi\colon U\to S$ restricted to each of the special fibres is smooth of relative dimension $1$ and that these restrictions induce finite morphisms on $Y\cap U$. Fibrewise criterion of flatness \cite[\href{https://stacks.math.columbia.edu/tag/039C}{Tag 039C}]{stacks-project} implies that $\pi\colon U\to S$ is flat, and the fibrewise criterion of smoothness \cite[\href{https://stacks.math.columbia.edu/tag/01V8}{Tag 01V8}]{stacks-project} implies that this morphism is smooth. By \cite[\href{https://stacks.math.columbia.edu/tag/0D4H}{Tag 0D4H}]{stacks-project}, the relative dimension of this morphism is $1$, because generically it is so, and the restriction of $\pi$ to the special fibres over $R$ is of relative dimension $1$. The openness of the quasi-finite locus \cite[\href{https://stacks.math.columbia.edu/tag/01TI}{Tag 01TI}]{stacks-project} implies that, possibly by shrinking $S$, we may assume that $\pi$ restricted to $Y\cap\pi^{-1}(S)$ is quasi-finite. To ensure that $Y\cap U=Y\cap\pi^{-1}(S)$, we spread out from the semilocalisation at $\pi(x_1),\ldots,\pi(x_n)$ to a small enough neighbourhood $S'\subset S$ such that $Y\cap U\supseteq Y\cap S'$, and then replace $S$ by this $S'$. Finally, Zariski's main theorem \cite[\href{https://stacks.math.columbia.edu/tag/03GU}{Tag 03GU}]{stacks-project} shows that $\pi$ restricted to $Y\cap U=Y\cap\pi^{-1}(S)$ is finite.
	 \par 
	 \es
	 %As a point of departure of the proof of \Cref{thm:A} (see \Cref{intro:sketch}), we establish \Cref{thm:B} from the introduction.~Presentation Lemma~\ref{prop:presentation_lemma} is inspired from \ccite{ces_grothendieck-serre}{Proposition~4.1}, where it was proved in the case of semilocal Dedekind rings. We note that the codimension $\geqslant 2$ hypothesis on $Y$ is to ensure the finiteness of the morphism $\pi$ when restricted to $Y\cap U$ (which, as we mentioned in the introduction, is important in the proof of \Cref{thm:A}).  
	 %\par We are now in a position to prove our presentation lemma from the introduction.
%	 \begin{pl}\label{prop:presentation_lemma}
%	 	For 
%	 	\bun
%	 	\item[$\circ$] a semilocal Pr\"ufer ring $R$ of Krull dimension at most $1$, 
%	 	\item[$\circ$] a smooth $R$-scheme $X$ fibrewise of pure relative dimension $d>0$,
%	 	\item[$\circ$] a closed subscheme $Z\subset X$ that is of codimension $\geqslant 2$,\footnote{This is equivalent to the condition that the $R$-generic fibres of $Z$ are of codimension $\geqslant 2$ while its other $R$-fibres are of codimension $\geqslant 1$ in that of $X$.}
%	 	\item[$\circ$] a closed subscheme $Z\subset Y\subset X$ that is $R$-fibrewise of codimension $\geqslant 1$, and
%	 	\item[$\circ$] points $x_1,\ldots,x_n\in X$;
%	 	\eun 
%	 	there are affine opens $x_1,\ldots,x_n\in U\subset X$ and $S\subset \bb{A}^{d-1}_{R}$ and a smooth $R$-morphism $\pi\colon U\to S$ of pure relative dimension $1$ such that $\pi|_{Z\cap U}$ is finite and $\pi|_{Y\cap U}$ is quasi-finite. 
%	 \end{pl}
	 \bs[Proof of Presentation~Lemma~\ref{prop:presentation_lemma}]
	 	Since the claim is Zariski-local around $x_1,\ldots,x_n\in X$, we may, without loss of generality, localise $R$ at the images of $x_1,\ldots,x_n$ and then ultimately spread out to assume that $R$ is, in fact, semilocal. Restricting to a connected component of $\spec(R)$, we may further assume that $\spec(R)$ is connected, and hence, integral (see Remark~\ref{rem:Prufer_ring}\eqref{rem:prufer_ring_pt:1}). If $R$ is of Krull dimension $0$, it is therefore a field (Remark~\ref{rem:Prufer_ring}\eqref{rem:prufer_ring_pt:2}), for which the statement follows from Gabber's presentation lemma \cite{ct_hoobler_kahn, hogadi_kulkarni_gabber_finite_field}. Thus, it suffices to assume that $R$ is of Krull dimension one.
	 	\par Choosing an embedding of $X$ into some $R$-affine space, let $j\colon X\hookrightarrow\overline{X}$ be the schematic image of the corresponding morphism from $X$ to the $R$-projective space. This constructed scheme $\overline{X}$ is $R$-flat, a fact that follows from \cite[\href{https://stacks.math.columbia.edu/tag/01RE}{Tag 01RE}]{stacks-project} and \ccite{bourbaki}{Chapter~I, \textsection2.4, Proposition~3(ii)}. Furthermore, thanks to \ccite{raynaud_gruson}{premi\`ere partie, corollaire~3.4.7}, $\overline{X}$ is of $R$-finite presentation. Hence, by \cite[\href{https://stacks.math.columbia.edu/tag/02FZ}{Tag 02FZ} and \href{https://stacks.math.columbia.edu/tag/0D4J}{Tag 0D4J}]{stacks-project}, it is also of $R$-fibrewise of pure dimension $d$. The flatness of $X$ and the constancy of fibrewise dimension ensures that the special $R$-fibres of $X$ are of codimension $1$ in X (\cite[\href{https://stacks.math.columbia.edu/tag/0D4H}{Tag 0D4H}, cf.~\href{https://stacks.math.columbia.edu/tag/054L}{Tag 054L}]{stacks-project}). By \cite[\href{https://stacks.math.columbia.edu/tag/081I}{Tag 081I}]{stacks-project}, the generic fibre of $X$ is dense in $\overline{X}$. This implies that the special fibres of $\overline{X}$ are of codimension $\geqslant 1$ in $\overline{X}$, and on the other hand, they are of codimension $\leqslant 1$ thanks to \cite[\href{https://stacks.math.columbia.edu/tag/0D4I}{Tag 0D4I}]{stacks-project}; showing that they are of codimension $1$. 
		 \par Let $C\subset\spec(R)$ be subscheme of closed points. We define $\overline{Y}$ to be the schematic closure of $Y'\colonequals Z\cup Y_C$ in $\overline{X}$. Given that the points of $Y'$ are of height $\geqslant 2$, the points of $\overline{Y}\setminus Y'$ are of height $\geqslant 3$. The generic fibre of $\overline{Y}\setminus Y'$ is of codimension $\geqslant 3$, and since the special fibres of $\overline{X}$ are of codimension $1$ in $\overline{X}$, the special fibres of $\overline{Y}\setminus Y'$ are of codimension $\geqslant 2$ in the corresponding special fibres of $\overline{X}$. We apply Theorem~\ref{prop:moving_lemma_general} to $(\overline{X}, X, \text{ points }x_1\ldots,x_n, \overline{Y})$, i.e., by inputting our $X$ as the $W$ of the proposition, our $\overline{X}$ as the $X$, and our $\overline{Y}$ as the $Y$. Hence, possibly by replacing $X$ with an affine open neighbourhood $U\subset X$ of $x_1,\ldots,x_n$, there exist an affine open $S\subset \bb{A}^{d-1}_{R}$ and a smooth $R$-morphism $\pi\colon X\to S$ of pure relative dimension $1$ such that $\pi|_Z$ and $\pi|_{Y_C}$ are finite, and in particular, $\pi|_{Y_C}$ is at least quasi-finite. To conclude our proof, it remains to use the openness of the quasi-finite locus \cite[\href{https://stacks.math.columbia.edu/tag/01TI}{Tag 01TI}]{stacks-project} which ensures that $\pi|_{Y}$ is also quasi-finite.   %\rb{Constant fibrewise dimension}
	 \es
%	 \br
%	 	We note that the preceding arguments required $X$ to be catenary, which is the content of \Cref{lem:valuation_rings_catenary}.
%	 \er
%%%%%%%%%%%%%%%%%%%%%%%%%%%%%%%%%%%%%%%%%%%%%%%%%%%%%%%%%%%%%%%%%%%%%%%%%%%%%%%%%%%%%%%%%%%%%%%%%%%%%%%%%%%%%%%
\section{Totally isotropic case of Conjecture~\ref{conj:gro_serre_valuation} for valuation rings of rank $1$.}\label{section:main_proof}
	%We recall the following from \ccite{phd-thesis}{Lemma~3.10}.
	%\bl[cf.~\ccite{mb_valuation_ring}{Théorème~A}]\label{lem:snail_generic_point_of_special_fibre}
	%	Given a valuation ring $V$, an integral domain $A$ that is an essentially smooth, faithfully flat $V$-algebra, a prime $\fr{p}\subset A$ corresponding to a generic point of the $V$-special fibre of $\spec A$, the morphism $V\to A_{\fr{p}}$ is an extension of valuation rings which induces an isomorphism at the level of value groups.
	%\el
	In this section, our goal is to prove Theorem~\ref{thm:GS-totally_isotropic}, closely following the strategy of the proof of \ccite{ces-fedorov}{Theorem~4.3}. Our result is obtained as a consequence of Theorem~\ref{thm:GS-totally_isotropic_body}, which serves as the main technical statement of this section. A key step in its demonstration is provided by Proposition~\ref{prop:2.6}.
	\par We begin by introducing the notion of totally isotropic reductive group schemes. Let $G$ be a reductive group scheme over a scheme $S$ and $G^{\mathrm{ad}}\colonequals G/Z$ be its adjoint quotient, where $Z\subseteq G$ is the centre. 
%	\end{montobo}
		\bdf[\ccite{ces_grothendieck-serre}{Definition~8.1}]\label{defn:totally_istropic_groups} 
		The group $G$ is called \textit{isotropic} if it contains a copy of $\bb{G}_{m,S}$ as a subgroup. It is called \textit{totally isotropic} if for any $s\in S$, every factor $G_i$ in the canonical decomposition \bud G^{\mathrm{ad}}_{\ca{O}_{s,S}}=\prod_{i=1}^nG_i,\eud in the sense of \ccite{sga3iii}{exposé~XXVI, corollaire~6.12}, is the Weil restriction from a connected finite \'etale cover $S_i\to\ca{O}_{S,s}$ of an isotropic, adjoint $S_i$-reductive group whose geometric fibres are simple. 
		\edf
	The class of totally isotropic groups contains, but is not limited to, the class of quasi-split groups; that is, reductive groups that contain a Borel subgroup (see \ccite{ces_torsors}{\textsection1.3.6}).
%	\br
%		Thanks to the uniqueness of the canonical decomposition as simple groups in , the above definition coincides with .
%	\er
	 %Unlike in the abelian case, we cannot just use presentation lemma to reduce to the noetherian case by some sort of Nisnevich descent. 
	 \par In the proof of Proposition~\ref{prop:2.6}, the following results play an important role. Proposition~\ref{prop:semilocal_Prufer_Grothendieck_serre} shows the semilocal Prüfer ring case of Conjecture~\ref{conj:gro_serre_valuation}. Whereas, Proposition~\ref{prop:shang_li} will allow us to reduce to the constant group case. The last one, Proposition~\ref{prop:purity_for_torsors}, which proves purity for $G$-torsors, will be used to extend torsors to any point of height $\leqslant 2$.
	 \bp[\ccite{ning_liu_constant}{Appendix A} and \ccite{phd-thesis}{Theorem~5.11}]\label{prop:semilocal_Prufer_Grothendieck_serre}
	 	Given a semilocal Prüfer ring $R$ and a reductive $R$-group scheme $G$, any generically trivial $G$-torsor over $R$ is trivial. 
	 \ep
	 \bs
	 	Passing to a connected component of $R$, we may assume that $R$ is even connected, and hence, integral (see Remark~\ref{rem:Prufer_ring}\eqref{rem:prufer_ring_pt:1}). The rest is the content of loc.~cit.
	 \es
	 The following result was proved by \citeauthor{shang_li_equivariant_compactification}.
	 \bp[\ccite{shang_li_equivariant_compactification}{Proposition 7.5}]\label{prop:shang_li}
	 	Let $X$ be a normal scheme, $Z\subseteq X$ be a finitely presented closed subscheme and $G$ and $G'$ be reductive $X$-group schemes with equal root datum over each geometric fibre. Then, any isomorphism $G_Z\cong G'_Z$ over $Z$ lifts to an isomorphism $G_{\tilde{X}}\cong G'_{\tilde{X}}$ over a finite étale cover $\tilde{X}\to X$ equipped with a section $Z\to\tilde{X}$.
	 \ep
	 \bs
	 	By a limit argument, we may reduce to the case when $X$ is even noetherian. In this, semilocalising and then spreading out, it reduces to loc.~cit.
	 \es
	 The following was proved by \citeauthor{ct-sansuc_extend_torsor_tori} in \ccite{ct-sansuc_extend_torsor_tori}{Theorem 6.13} in the case when the base is a regular scheme.
	 \bp\label{prop:purity_for_torsors}
	 	Let $R$ be a Pr\"ufer domain, let $X$ be an ind-smooth, integral $R$-scheme and let $j\colon U\hookrightarrow X$ be a quasi-compact open such that at each point $z\in Z\colonequals X\setminus U$ lying over $y$, we have \bd\dim(\ca{O}_{X_{y},z})+\min(1,\dim(R_y))= 2.\ed
	 	Then, for any point $z\in Z$ and reductive $X$-group scheme $G$, the restriction morphism induces an equivalence
	 	\bud\label{diag:lifting_torsors}\mathbf{B}G(\spec(\ca{O}_{X,z}))\isom&\mathbf{B}G(\spec(\ca{O}_{X,z})\setminus\{z\}).\eud
	 \ep
	  The case when $G$ is a torus was established in \ccite{phd-thesis}{Proposition~4.5} (see also \ccite{amar:unramified_case_gersten_conjecture}{Proposition~4.3}). Whereas, for an arbitrary reductive $G$, the result was proved independently in \cite{ning_liu_quasi-split}.
	 \bs[Sketch of proof of Proposition~\ref{prop:purity_for_torsors}]
	 	As in \cite{ct-sansuc_extend_torsor_tori}, through a faithful embedding $G\hookrightarrow\GL_n$ whose cokernel is necessarily affine, we may first reduce to the case when $G=\GL_n$. In this case, one can use \ccite{phd-thesis}{Lemma~4.4} (see also \ccite{amar:unramified_case_gersten_conjecture}{Lemma 4.2}) to show that the extension of a locally free sheaf is reflexive, which, thanks to \ccite{gabber-ramero-foundations}{Proposition~11.4.1(iii)}, must even be locally free.
	 \es
	The final ingredient in Proposition~\ref{prop:2.6} is the following result, due to \citeauthor{ces-fedorov} in \cite{ces-fedorov}, which demonstrates the generalised Horrocks' principle for torsors under a totally isotropic group over any ring.
	\bp[\ccite{ces-fedorov}{Theorem~4.2}]\label{lem:generalised_horrocks_principle}
	Given a ring $A$ and a totally isotropic reductive $A$-group scheme $G$, any $G$-torsor $\scr{E}$ over $\bb{P}^1_A$ such that $\scr{E}|_{\{t=\infty\}}$ is trivial actually trivialises over $\bb{A}^1_A$.
	\ep
	 \par We are now in a position to establish a key step in the proof of Theorem~\ref{thm:GS-totally_isotropic_body}, following the arguments in \ccite{ces-fedorov}{Proposition~2.6}.
	 \bp\label{prop:2.6}
	 	Let 
	 	\begin{itemize} 
	 			\item[$\circ$] $R$ be a Prüfer ring,
	 			\item[$\circ$] $X$ be a smooth $R$-scheme with points $x_1,\ldots,x_n\in X$,
	 			\item[$\circ$] $G$ be a totally isotropic reductive $X$-group scheme, and
	 			\item[$\circ$] $E$ be a generically trivial $G$-torsor over $X$.
	 	\end{itemize}
	 	Then, Zariski semilocally around $x_1,\ldots,x_n\in X$, there exist
	 	\begin{enumerate}[(i)]
	 		\item\label{prop:2.6:part:1} principal closed subscheme $Y\subset X$ that is $R$-fibrewise of positive codimension,
	 		\item\label{prop:2.6:part:2} closed subscheme $Z\subset Y$ such that the $R$-generic fibres of $Z$ are even of codimension $\geqslant 2$ in those of $X$, and %\textcolor{gray}{an open $C\subseteq\bb{P}^1_X$ containing $\bb{P}^1_{X\setminus Z}$ along with the sections $\{t=0\}$ and $\{t=\infty\}$ in $\bb{P}^1_X(X)$, and}
	 		\item\label{prop:2.6:part:3} a $G$-torsor $\scr{E}$ over $\bb{P}^1_{X\setminus Z}$ such that $\scr{E}|_{\{t=0\}}\cong E|_{X\setminus Z}$ whereas $\scr{E}|_{\{t=\infty\}}$ and $\scr{E}|_{\bb{P}^1_{X\setminus Y}}$ are trivial.
	 	\end{enumerate}
	 \ep
	 \bs
	 	Since the claim is Zariski-semilocal around $x_1,\ldots x_n\in X$, we may, without loss of generality, assume that $R$ is semilocal by semilocalising at the images of these points. For the same reason, by restricting to a connected component of $\spec(R)$ (resp., of $X$), we may even assume that $\spec(R)$ (resp., $X$) is connected. Consequently, by Remark~\ref{rem:Prufer_ring}\eqref{rem:prufer_ring_pt:1}, $R$ is a Prüfer domain and, since $X$ is $R$-smooth, it is integral. 
	 	\par Moreover, if $E$ is trivial, we may take $C=\bb{P}^1_X$ and let $\scr{E}$ be the trivial $G$-torsor, thereby proving the claim. If the relative dimension $d$ of $X$ over $R$ equals $0$, then $X$ itself is the spectrum of a semilocal Prüfer domain. In this situation, Proposition~\ref{prop:semilocal_Prufer_Grothendieck_serre} shows that $E$ is trivial, therefore settling the claim. We are thus reduced to the case $d\geqslant 1$. 
	 	\par In addition, by a limit argument, we may further restrict to the case when $R$ is of finite Krull dimension. Indeed, \ccite{amar:unramified_case_gersten_conjecture}{Lemma~2.5(b)} proves that $R$ is an increasing union of its Prüfer sub-domains $R_{\lambda}$ of finite Krull dimension. These canonical maps $R_{\lambda}\hookrightarrow R$ are, in fact, flat (see Remark~\ref{rem:Prufer_ring}\eqref{rem:prufer_ring_pt:3}), and since flat morphisms preserve fibrewise dimension, we may descend all data to some $R_{\lambda}$, and then base change back to $R$. %we reduce to proving the claim for $P$ of finite Krull dimension.
	 	\par Since the semilocalisation $\ca{O}$ of $X$ at the set of generic points of the $R$-special fibres of $X$ is a semilocal Prüfer domain (see \ccite{amar:unramified_case_gersten_conjecture}{Lemma 2.13}), again, by Proposition~\ref{prop:semilocal_Prufer_Grothendieck_serre}, the generically trivial $G$-torsor $E|_{\spec(\ca{O})}$ must be trivial. By a limit argument, part (2) of loc.~cit.~then allows us to find an element $f\in\Gamma(X,\ca{O}_{X})$ that maps to a unit in $\ca{O}$, whose vanishing locus $H_f\subseteq X$ is $R$-flat and at the same time, our $E$ restricts to a trivial $G$-torsor over $X\setminus H_f$. In particular, this scheme $H_f$ is even $R$-fibrewise of positive codimension in $X$. Thanks to \ccite{amar:unramified_case_gersten_conjecture}{Presentation Lemma~3.2}, Zariski-semilocally around $x_1,\ldots,x_n\in X$, there exist an affine open $S\subseteq\bb{A}^{d-1}_R$ and a smooth $R$-morphism $\pi\colon X\to S$ of relative dimension $1$ such that $\pi|_{H_f}$ is quasi-finite. Let $A$ be the semilocalisation of $X$ at $x_1,\ldots,x_n$ and $\lambda\colon\spec(A)\to X$ be the canonical morphism. Base changing $\pi$ along $\pi\circ\lambda$, we obtain a smooth $A$-scheme $C$ of pure dimension $1$ along with an $A$-quasi-finite closed subscheme $H\subseteq C$ and a section $s\in C(A)$ lifting $\lambda$. Furthermore, we obtain a reductive $C$-group scheme $\scr{G}$ satisfying $s^{\ast}\scr{G}\cong G$ and a $\scr{G}$-torsor $\scr{E}$ such that $s^{\ast}\scr{E}\cong E$. 
	 	A priori, this $C$ need not even embed inside $\bb{P}^1_A$ and $\scr{G}$ need not equal the base change $G_C$ of the $A$-scheme $G$. We address these issues below.
	 	\item Passing to a finite étale cover of $C$, by Proposition~\ref{prop:shang_li}, we can immediately reduce to the case when $\scr{G}=G_C$. Then, over each maximal ideal $\fr{m}\subseteq A$, by removing a principal closed subscheme containing the points in $H_{\kappa(\fr{m})}\setminus s(\kappa(\fr{m}))$ (\stacks{00ds}), where we note that $H_{\kappa(\fr{m})}$ is a finite discrete set of points since it is quasi-finite over a field, we produce an affine open $C'\subset C$ containing $s(A)$ such that $(H\cap C')_{\kappa(\fr{m})}= s(\kappa(\fr{m}))$ over each maximal ideal $\fr{m}\subseteq A$. Consequently, $H'\colonequals H\cap C'\subseteq C'$ has no \textit{finite field obstruction} to embedding it in any projective space over $A$, in the sense of \ccite{ces-fedorov}{Definition~2.3}, since the same holds for $\spec(A)$ itself. Therefore, thanks to \ccite{ces-fedorov}{Lemma 2.5}, there are an affine open $C''\subseteq C'\sqcup\bb{A}^1_A$ containing $s(A)\sqcup\{t=0\}$, an affine open $\tilde{C}\subseteq\bb{A}^1_A$ and an étale morphism $f\colon C''\to\tilde{C}$ that maps $H\sqcup\{t=0\}$ isomorphically onto a closed subscheme $\tilde{H}\subseteq\tilde{C}$ such that $$(H\sqcup\{t=0\})=\tilde{H}\times_{\tilde{C}}C''.$$ This simultaneously ensures that $s\in C'(A)$ descends to a section $\tilde{s}\in\tilde{C}(A)$ and $\{t=0\}\in\bb{A}^1_A$ descends to a section $\tilde{s}_0\in\tilde{C}(A)$ disjoint from $\tilde{s}$. 
	 	Thus, extending $\scr{E}|_{C'}$ to a $G$-torsor $\scr{E}''$ over $C''$ such that $\scr{E}''|_{\bb{A}^1_A\cap C''}$ is trivial and then patching $\scr{E}''$ with a trivial $G$-torsor over $\tilde{C}\setminus\tilde{H}$, we obtain a $G$-torsor $\tilde{\scr{E}}$ over $\tilde{C}$ satisfying that $\tilde{s}^{\ast}(\tilde{\scr{E}})=E$ as well as $\tilde{s}^{\ast}_0\tilde{\scr{E}}$ and $\tilde{\scr{E}}|_{\tilde{C}\setminus\tilde{H}}$ are trivial. 
	 	Zooming in the generic fibre, \ccite{gille_torseurs_sur_la_droite_affine}{corollaire 3.10(a)} ensures further that $\tilde{\scr{E}}|_{\tilde{C}_{K}}$ is trivial, where $K$ is the fraction field of $A$. 
	 	\par Hence, we may replace $C$ by this $\tilde{C}$, and their corresponding data, and finally reduce to the case when $C\subseteq\bb{A}^1_A$. 
	 	 At this point, however, we have two disjoint sections in $C(A)$ at our disposal. By modifying the coordinates of $\bb{A}^1_A$ we arrange further that $s$ (resp., $s_0$) is $\{t=0\}$ (resp., $\{t=\infty\}$). Indeed, first, acting by an automorphism $g$ of $\bb{P}^1_A$, which comes at the cost of replacing $C\subset\bb{P}^1_A$ and all its dependants by their respective images under $g$, we may transport $s_0$ to the section $\{t=\infty\}$. Then, using a suitable automorphism of $\bb{P}^1_A$ which fixes $\infty$, which amounts to a linear change of coordinates, we may additionally move $s$ to the section $\{t=0\}$. %thereby establishing the claims \eqref{prop:2.6:part:2} and \eqref{prop:2.6:part:3} for the semilocal ring $R$, instead of $X$, as required. Below we prove the remaining claim, i.e., claim \eqref{prop:2.6:part:1}, over $R$.
	 	\par To prove the claims \eqref{prop:2.6:part:1}-\eqref{prop:2.6:part:3}, we shall prove their counterparts over $A$ and then spread them out to a Zariski-semilocal affine neighbourhood of $x_1,\ldots,x_n\in X$.  Since $\scr{E}$ trivialises over $C_{K}$, we can assume that it does so also over $C_T$, where $T\subseteq\spec(A)$ is a dense open subset. Extending $\scr{E}$ to $C\cup\bb{P}^1_T$ by patching it with the trivial $G$-torsor over $\bb{P}^1_T$, we may assume, without loss of generality, that $\bb{P}^1_T\subseteq C$. 
	 	In particular, since $C$ is $A$-fibrewise dense in $\bb{P}^1_A$, it follows that $H$ is $A$-fibrewise dense in its closure $\overline{H}$ in $\bb{P}^1_A$, showing that $\overline{H}$ is $A$-finite. Again, patching $\scr{E}$ with a trivial torsor over $\bb{P}^1_A\setminus\overline{H}$, we may assume that $(\bb{P}^1_A\setminus\overline{H})\subseteq C$. Therefore, the complement $B\colonequals\bb{P}^1_A\setminus C$ does not contain any irreducible component of an $A$-fibre of $\bb{P}^1_A$, showing that $B$ neither contains any irreducible component of an $R$-fibre of $\bb{P}^1_A$. Additionally, since also $\bb{P}^1_K\subseteq C$, any height $\leqslant 1$ point in the $R$-generic fibre of  $\bb{P}^1_A\setminus\bb{P}^1_{K}$, which must be the generic point of some $A$-fibre of $\bb{P}^1_A$, is also contained in $C$. 
	 	Proposition~\ref{prop:purity_for_torsors} then allows us extend $\scr{E}$ to any height $\leqslant 2$ point in the $R$-generic fibre of $\bb{P}^1_{A}$ and to any $\leqslant 1$ point in other $R$-fibres of $\bb{P}^1_A$. Consequently, we may further assume that the $R$-fibres of $B$ are of codimension $\geqslant 2$ in $\bb{P}^1_A$, in addition, its $R$-generic fibre is even of codimension $\geqslant 3$ in $\bb{P}^1_A$. 
	 	\par As a consequence, the dimension formula \stacks{02JU}, applied $R$-fibrewise, guarantees that the image $Z_A\subset\spec(A)$ of $B$ is a closed subscheme that is $R$-fibrewise of codimension $\geqslant 1$, meanwhile its $R$-generic fibre is even of codimension $\geqslant 2$. 
	 	At the same time, by definition, our $C$ actually contains $\bb{P}^1_{\spec(A)\setminus Z_A}$. Thanks again to \ccite{amar:unramified_case_gersten_conjecture}{Lemma 2.13(2)}, we find a principal closed subscheme $Y_A\subset\spec(A)$ containing $Z_A$ that is $R$-fibrewise of positive codimension (in fact, it is even flat over $R$). This proves \eqref{prop:2.6:part:1}-\eqref{prop:2.6:part:2} over $A$. 
	 	\par Replacing $C$ by $\bb{P}^1_{\spec(A)\setminus Z_A}$, by construction, we have $\scr{E}|_{\{t=0\}}\cong E|_{\spec(A)\setminus Z_A}$ and $\scr{E}|_{\{t=\infty\}}$ is trivial. To prove \eqref{prop:2.6:part:3} over $A$, it suffices to demonstrate that $\scr{E}|_{\bb{P}^1_{\spec(A)\setminus Y_A}}$ is trivial. However, a priori, this might not hold. We shall replace $Z_A$ by a larger subscheme $Z'_A\subset Y_A$, as well as modify our $\scr{E}$ to resolve this issue.
	 	\par Since $\spec(A)\setminus Y_A$ is affine and $\scr{E}|_{\{t=\infty\}}$ is trivial, Proposition~\ref{lem:generalised_horrocks_principle} shows that $\scr{E}|_{\bb{A}^1_{\spec(A)\setminus Y_A}}$ is trivial, whence, again by Proposition~
	 	\ref{lem:generalised_horrocks_principle}, it follows that $\scr{E}|_{\bb{P}^1_{\spec(A)\setminus Y_A}\setminus\{t=1\}}$ is trivial. Thus, by patching $\scr{E}$ with a trivial $G$-torsor over $C\setminus\{t=1\}$, we obtain a $G$-torsor $\scr{E}'$ over $C'\colonequals (C\setminus\{t=1\})\cup\bb{P}^1_{\spec(A)\setminus Y_A}$. Considering $B'\colonequals\bb{P}^1_A\setminus C'$, in a similar vein as above, by Proposition~\ref{prop:purity_for_torsors}, there exists a closed subscheme $Z_A\subseteq Z'_A\subset Y_A$ such that $\scr{E}'$ extends all the way to $\bb{P}^1_{\spec(A)\setminus Z'_A}$. Consequently, replacing $Z_A$ by $Z'_A$ and $\scr{E}$ by this extension of $\scr{E}'$, we may further assume that $\scr{E}|_{\bb{P}^1_{\spec(A)\setminus Y_A}}$ is trivial, establishing \eqref{prop:2.6:part:3} as well. Thus, we are done.
	 	%\par Therefore, to prove \eqref{prop:2.6:part:1}-\eqref{prop:2.6:part:2}, it suffices to spread $Z_R\subset Y_R\subset\spec(R)$ in a Zariski-semilocal affine neighbourhood of $x_1,\ldots,x_n\in X$ to $Z\subset Y\subset X$, as required. To prove the claim over $X$, it remains for us to spread out the objects defined over $R$ to their counterparts over a Zariski neighbourhood of $x_1,\ldots x_n\in X$. Again, by abuse of notation, likewise, we produce an open $C\subseteq\bb{P}^1_X$ containing $\bb{P}^1_{X\setminus Z}$, for some closed $Z\subset X$ whose $P$-generic fibre is of codimension $\geqslant 2$ while other $P$-fibres are of codimension $\geqslant 1$, together with disjoint sections $s,s_0\in\bb{P}^1_X(X)$ and a $G$-torsor $\scr{E}$ over $C$ such that $s^{\ast}\scr{E}\cong E$ and $s_0^{\ast}\scr{E}$ as well as $E|_{\bb{P}^1_{K}}$ are trivial. \rl{Need to produce $Y$ now.}
	 \es
%		The following is the generalised Horrocks' principle for torsors under a totally isotropic group over any ring.
%	 \bp[\ccite{ces-fedorov}{Theorem~4.2}]\label{lem:generalised_horrocks_principle}
%	 	Given a ring $A$ and a totally isotropic reductive $A$-group scheme $G$, any $G$-torsor $\scr{E}$ over $\bb{P}^1_A$ such that $\scr{E}|_{\{t=\infty\}}$ is trivial actually trivialises over $\bb{A}^1_A$.
%	 \ep
	 Let us prove the main technical result of this section. Our proof is based on \ccite{ces-fedorov}{Theorem~4.3}.
	 \bt\label{thm:GS-totally_isotropic_body}
	 Let $R$ be a Prüfer ring of Krull dimension $1$, let $A$ be the semilocalisation of a smooth $R$-algebra at finitely many points, and let $G$ be a totally isotropic reductive $A$-group scheme. Then, given a generically trivial $G$-torsor $E$ over $A$, there exist 
	 \begin{enumerate}[(i)]
	 	\item an open $C\subseteq\bb{A}^1_A$ with a section $s\in C(A)$,
	 	\item a $G$-torsor $\scr{E}$ over $C$ such that $s^{\ast}\scr{E}\cong E$, 
	 	\item an $A$-finite closed subscheme $\scr{Z}\subset C$, and
	 	\item  a $G$-torsor $\tilde{\scr{E}}$ over $\bb{P}^1_{C\setminus\scr{Z}}$ such that $\tilde{\scr{E}}|_{\{t=0\}}\cong\scr{E}|_{C\setminus\scr{Z}}$ while $\tilde{\scr{E}}|_{\{t=\infty\}}$ is trivial.
	 \end{enumerate}
	 \et
	 \bs
	 	 We suppose that $X$ is a smooth $R$-scheme such that $A$ is the semilocalisation of $X$ at points $x_1,\ldots,x_n\in X$. Possibly by shrinking $X$ around $x_1,\ldots,x_n$, we may also assume that $G$ as well as $E$ begin life over $X$ itself. Furthermore, restricting ourselves to a connected component of $A$ (resp., of $X$), we may suppose that $A$ (resp., $X$) is connected. Since $X$ is $R$-smooth, it follows then that $X$, and hence, $A$ is integral. Moreover, since the claim is Zariski-local around $x_1,\ldots,x_n\in X$, we can semilocalise $R$ at the images of $x_1,\ldots,x_n$ to restrict ourselves to the case when $R$ is semilocal.
	 	 %\par To prove the claim, it suffices to show that $E$ is trivial in a Zariski neighbourhood of $x_1,\ldots,x_n\in X$. In particular, zooming in a connected component, we may assume that $X$ and hence $R$ is connected. Since both $X$ and $R$ reduced, they are even integral. In the same spirit, we may localise $P$ at the images of $x_1,\ldots,x_n$ to arrange that is is semilocal and also integral. 
	 	 \par Proposition~\ref{prop:2.6} then furnishes us with closed subschemes $Z\subset Y\subset X$ such that $Y$ is a principal closed subscheme that is $R$-fibrewise of positive codimension in $X$, and additionally, the $R$-generic fibres of $Z$ are of codimension $\geqslant 2$ in those of $X$, as well as $$\text{a $G$-torsor $\scr{E}_0$ over $\bb{P}^1_{X\setminus Z}$ such that $\scr{E}_0|_{\{t=0\}}\cong E|_{X\setminus Z}$ while $\scr{E}_0|_{\{t=\infty\}}$ and $\scr{E}_0|_{\bb{P}^1_{X\setminus Y}}$ are trivial.}$$ 
	 	 %an open $F\subseteq\bb{P}^1_X$ containing $\bb{P}^1_{X\setminus Z}$ as well as sections $\{t=0\}$ and $\{t=\infty\}$ in $\bb{P}^1_X(X)$ and a $G$-torsor $\scr{E}$ over $F$ such that $\scr{E}|_{\{t=0\}}\cong E$ while $\scr{E}|_{\{t=\infty\}}$ and $\scr{E}|_{\bb{P}^1_{X\setminus Y}}$ are trivial. In particular, if $F$ were equal to $\bb{P}^1_X$, \Cref{lem:section_invariance} would already demonstate that $E$ is trivial, but, a priori, this might not be the case. Henceforth, we arrange our data so that ultimately \Cref{lem:section_invariance} may be used.
	 	 \par  Letting $d$ be the relative dimension of $X$ over $R$, we note that Proposition~\ref{prop:semilocal_Prufer_Grothendieck_serre} establishes the claim when $d=0$. Therefore, it suffices to treat the case $d\geqslant 1$. Consequently, thanks to Presentation~Lemma~\ref{prop:presentation_lemma}, Zariski-semilocally around $x_1,\ldots x_n\in X$, there are an affine open $S\subseteq\bb{A}^{d-1}_R$ and a smooth morphism $\pi\colon X\to S$ of pure relative dimension $1$ such that $\pi|_Y$ and $\pi|_Z$ are $S$-quasi-finite and $S$-finite respectively. 
	 	 \par Let $\lambda\colon\spec(A)\to X$ be the canonical morphism induced by semilocalisation at $x_1,\ldots,x_n\in X$. Base changing $\pi$ along $\pi\circ\lambda$, we obtain a smooth $A$-scheme $C$ of pure dimension $1$ along with an $A$-quasi-finite closed subscheme $\scr{Y}\subseteq C$, an $A$-finite closed subscheme $\scr{Z}\subseteq C$ and a section $s\in C(A)$ lifting $\lambda$. Additionally, pulling-back the $X$-group scheme $G$, we obtain a reductive $C$-group scheme $\scr{G}$ satisfying $s^{\ast}\scr{G}=G$ and a $\scr{G}$-torsor $\scr{E}$ over $C$ such that $s^{\ast}\scr{E}=E$. On top of that, base changing the $G$-torsor $\scr{E}_0$ over $\bb{P}^1_{X\setminus Z}$, we get a $\scr{G}$-torsor $\tilde{\scr{E}}$ over $\bb{P}^1_{C\setminus\scr{Z}}$ such that $\tilde{\scr{E}}|_{\{t=0\}}\cong\scr{E}|_{C\setminus\scr{Z}}$ whereas $\tilde{\scr{E}}|_{\{t=\infty\}}$ as well as $\tilde{\scr{E}}|_{\bb{P}^1_{C\setminus\scr{Y}}}$ are trivial.
	 	 \par Our goal now is to arrange $C$ to be an open of $\bb{A}^1_A$. In a similar vein as the proof of Proposition~\ref{prop:2.6}, replacing $C$ by open containing $\scr{Y}\cup s(A)$, we first reduce to the case when $\scr{G}=G_C$ (see Proposition~\ref{prop:shang_li}). Thereafter, passing to a finite étale cover of $C$ with a lift of the section $s$, which is also denoted abusively by $(C,s)$, and then an open of $C$ containing the pullback of $\scr{Y}$ and $s(A)$,  we can ensure that there is no finite field obstruction to embedding $\scr{Y}\cup s$ itself into $\bb{A}^1_A$. Consequently, possibly by replacing $C$ with an affine open containing $s(A)$, thanks to \ccite{ces-fedorov}{Lemma 2.5}, there are an affine open $C'\subseteq\bb{A}^1_A$ and an étale morphism $f\colon C\to C'$ that maps $\scr{Y}\cup s(A)$ isomorphically onto a closed subscheme $\scr{Y}'\subseteq C'$ such that $\scr{Y}\cup s(A)=\scr{Y}'\times_{C'}C$. Let $s'\in C'(A)$ be the section induced by $s$ and $\scr{Z}'\subset C'$ be the isomorphic image of $\scr{Z}$. Thus, patching $\scr{E}$ with a trivial $G$-torsor over $C'\setminus\scr{Y}'$ produces a $G$-torsor $\scr{E}'$ over $C'$ that trivialises away from $\scr{Y}'$ such that $s'^{\ast}\scr{E}'=E$. 
	 	 \par However, in order to descend $\tilde{\scr{E}}$ to a $G$-torsor $\tilde{\scr{E}}'$ over $\bb{P}^1_{C'\setminus\scr{Z}'}$, we first need to ensure that the trivialisation of $(\tilde{\scr{E}}|_{\{t=\infty\}})|_{C\setminus\scr{Y}}$ extends to a trivialisation of $\tilde{\scr{E}}|_{\bb{P}^1_{C\setminus\scr{Y}}}$. For instance, this can be done provided the restriction to $\{t=\infty\}$ induces a bijection on the set of all trivialisations, i.e., $\{t=\infty\}^{\ast}\colon G(\bb{P}^1_{C\setminus\scr{Y}})\iso G(C\setminus\scr{Y})$. Since $C\setminus\scr{Y}$ is affine, this follows from the fully faithfulness in \ccite{ces-fedorov}{Proposition~3.1(a)}. 
	 	 \par Thus, patching $\tilde{\scr{E}}$ with a trivial $G$-torsor over $\bb{P}^1_{C'\setminus\scr{Y}'}$, we obtain a $G$-torsor $\tilde{\scr{E}}'$ over $\bb{P}^1_{C'\setminus\scr{Z}'}$ that trivialises on $\bb{P}^1_{C'\setminus\scr{Y}'}$ such that $\tilde{\scr{E}}'|_{\{t=\infty\}}$ is trivial. However, the $G$-torsor $\scr{E}''\colonequals\tilde{\scr{E}}'|_{\{t=0\}}$ over $C'\setminus\scr{Z}'$ might actually be different from $\scr{E}'|_{C'\setminus\scr{Z}'}$ even though its pullback to $C\setminus\scr{Z}$ is certainly $\scr{E}|_{C\setminus\scr{Z}}$. Consequently, to reduce to the case when $C=C'$, particularly, when $C$ is an open of $\bb{A}^1_A$, it suffices to demonstrate that $\scr{E}''$, which, a priori, lives only over $C'\setminus\scr{Z}'$, actually extends all the way to a $G$-torsor over $C'$, in which case, we might as well replace $\scr{E}'$ by this $\scr{E}''$. 
	 	 \par To accomplish this, we patch $\scr{E}''$ with $\scr{E}$ along $C\setminus\scr{Z}$ to obtain a $G$-torsor $\scr{E}'$, denoted abusively, over all of $C'$, as required. Thus, we are done.
	 \es
	 Before demonstrating Theorem~\ref{thm:GS-totally_isotropic}, we recall the following result, also due to \citeauthor{ces-fedorov}, proven in \cite{ces-fedorov}, which establishes the sectional invariance of torsors over $\bb{P}^1_A$ over any semilocal ring $A$. This result plays a crucial role in the remainder of the proof.
	 \bp[\ccite{ces-fedorov}{Theorem~3.6}]\label{lem:section_invariance}
	 Given a semilocal ring $A$, a reductive $A$-group scheme $G$, a $G$-torsor $\scr{E}$ over $\bb{P}^1_A$; for any two sections $s_1,s_2\in\bb{P}^1_A(A)$, there is an isomorphism $s_1^{\ast}\scr{E}\simeq s_2^{\ast}\scr{E}$.
	 \ep
	 \bs[Proof of Theorem~\ref{thm:GS-totally_isotropic}]
	 	Let $E$ be a generically trivial $G$-torsor over $A$ that we need to prove is trivial. Thanks to Theorem~\ref{thm:GS-totally_isotropic_body}, whose notations we adopt here, we get an open $C\subseteq\bb{A}^1_A$ with a section $s\in C(A)$, which ensures, in particular, that $\bb{P}^1_A\setminus C$ is $A$-finite. Thus, the avoidance lemma \ccite{gabber-liu-lorenzini}{Theorem~5.1} provides a hyperplane $\scr{Z}\subset H\subset C$ that is closed in $\bb{P}^1_A$, and at the same time, not containing any generic point of an $A$-fibre of $C$. As a consequence, this $H$ is $A$-finite. 
	 	\par On the other hand, since $C\setminus H$ is affine, Proposition~\ref{lem:generalised_horrocks_principle} demonstrates that $\tilde{\scr{E}}|_{\bb{A}^1_{C\setminus H}}$ is trivial, as the same is true for $\tilde{\scr{E}}|_{\{t=\infty\}}$. Therefore, triviality is true even for $\scr{E}|_{C\setminus H}\cong(\tilde{\scr{E}}|_{\bb{A}^1_{C\setminus H}})|_{\{t=0\}}$. 
	 	\par Thereafter, patching $\scr{E}$ with a trivial $G$-torsor over $\bb{P}^1_A\setminus H$, we extend it to a $G$-torsor over the entire $\bb{P}^1_A$, while simultaneously ensuring that $\scr{E}|_{\{t=\infty\}}$ is trivial. In this case, Proposition~\ref{lem:section_invariance} shows that $E\cong s^{\ast}\scr{E}$ is trivial as well, establishing the claim.
	 	%Let $E$ be a generically trivial $G$-torsor over $R$ that we need to prove is trivial. Thanks to \Cref{thm:GS-totally_isotropic_body}, whose notations we adopt here, we get an open $C\subseteq\bb{A}^1_R$ with a section $s\in C(R)$, which ensures that $\bb{P}^1_R\setminus C$ is $R$-finite. Thus, the avoidance lemma \ccite{gabber-liu-lorenzini}{Theorem~5.1} provides a hyperplane $\scr{Z}\subset H\subset\bb{C}$ closed in $\bb{P}^1_R$, and at the same time, not containing any generic point of an $R$-fibre of $C$; particularly, this $H$ is $R$-finite. Since $C\setminus H$ is affine, the generalised Horrock's principle \ref{lem:generalised_horrocks_principle} demonstrates that $\tilde{\scr{E}}|_{\bb{A}^1_{C\setminus H}}$ is trivial because the same is true for $\tilde{\scr{E}}|_{\{t=\infty\}}$. In particular, triviality holds even for $\scr{E}|_{C\setminus H}\cong(\tilde{\scr{E}}|_{\bb{A}^1_{C\setminus H}})|_{\{t=0\}}$. Thereafter, patching $\scr{E}|_{C\setminus H}$ with a trivial $G$-torsor over an open neighbourhood of $H\cup\{t=\infty\}\subseteq\bb{P}^1_R$, we extend it to $G$-torsor over the entire $\bb{P}^1_R$ while ensuring that it satisfies $\scr{E}|_{\{t=\infty\}}$. In this case, the sectional invariance \ref{lem:section_invariance} shows that $E\cong s^{\ast}\scr{E}$ is trivial as well, establishing the claim.
	 \es
\printbibliography

\end{document}

%% file: preamble.tex
% !TeX spellcheck = en_GB
%\documentclass[11pt,a4paper,onesided, notitlepage]{article}
%\usepackage{fullpage} %To set narrower margins automatically
%\usepackage[tmargin=2.5cm,bmargin=3cm,lmargin=1.5cm,rmargin=1.25cm]{geometry}%Margins I use for research statements
\usepackage[tmargin=2cm,bmargin=2.5cm,lmargin=1.25cm,rmargin=1.25cm]{geometry}%Margins I use for the articles
\usepackage[utf8]{inputenc}
\usepackage[english]{babel}
\usepackage[rm]{titlesec}%Section title unbolden format. Remove when providing referee reports
\titlelabel{\thetitle.\quad}%Put a fullstop after the number preceeding the name of a section 1. 2. Remove when providing referee reports
%\renewcommand{\contentsname}{\hspace{2cm} Table of Contents\hfill}%Centering name of table of contents. Didn't quite work. 
%\renewcommand{\cftaftertoctitle}{\hfill}
%----------Packages-----------------------------------------------------
\usepackage[T1]{fontenc} %This is to include accented characters like ö and é
\usepackage{amssymb,amsfonts}
\usepackage{amsthm}
\usepackage{tikz-cd}
\usepackage[colorinlistoftodos]{todonotes}
\usepackage[framemethod=TikZ]{mdframed} %রঙিন Theorem ENv। মনে হয় amsthm এর আগে লোড করতে হবে।
\usepackage{calligra,mathrsfs}%---For mathscr
\usepackage{amsmath}
\usepackage{mathtools}%-------Loads amsmath
\usetikzlibrary{shapes.geometric, arrows}%---libraries loaded in tikz
\usepackage[backend=biber,style=alphabetic,sorting=nyt, maxnames=10, maxalphanames=10]{biblatex}%----biber(recommended) or bibtex(traditional) is the backbender that sorts the bib items, common styles are "numeric" and "alphabetic" and "draft", common sorting counter are "none" for no sorting
%Reduce the font of the references
\usepackage{csquotes}%When using babel or polyglossia with biblatex, loading csquotes is recommended to ensure that quoted texts are typeset according to the rules of your main language.
\usepackage{colonequals}%To get \colonequals instead of \coloneqq
\usepackage{adjustbox} %Don't remember anymore
\usepackage[colorlinks, final]{hyperref} %This is for hyperlink, links to references within the document. use citebordercolor, linkbordercolor, urlbordercolor for borders around links. For no borders put color white
\usepackage{cleveref}  %For \cref and \Cref
\usepackage{enumerate} %to get enumerate lists with roman and small roman letters
%\usepackage{enumitem}%to get cuter unordered lists
%%%%%%%%%%%%%%%%%%%%Ning's Colours%%%%%%%%%%%%%%%%%%%%%
\definecolor{imperialred}{RGB}{237, 41, 57}
\definecolor{royalblue}{RGB}{64, 106, 212}
\definecolor{link}{RGB}{11,0,128}
\hypersetup{
	colorlinks = true,
	breaklinks = true,
	citecolor=royalblue,
	linkcolor=imperialred,
	urlcolor=link,
	linktocpage = true
}
%\hypersetup{
%  colorlinks = true,
%  breaklinks = true,
%  linkcolor  = green, %ForestGreen,
%  citecolor  = blue, %ProcessBlue,
%  urlcolor   = purple, %RoyalPurple,
%  linktocpage = true,
%}     
%--------Theorem Environments--------------------------------------------
\makeatletter

\makeatother

\newrobustcmd*{\citeauthorfullname}{\AtNextCite{\DeclareNameAlias{labelname}{given-family}}\citeauthor}

%%%%%%%%%%%%%%%%%%%%%%%%%%%%%%%%%%%%%%%
\newtheorem{thm}{Theorem}[section]
\newtheorem{cor}[thm]{Corollary}
\newtheorem{prop}[thm]{Proposition}
%\newmdtheoremenv{statement}{Statement}[section]
%\newmdtheoremenv{lem}[thm]{Lemma}
\newtheorem{lem}[thm]{Lemma}
\newtheorem{conj}[thm]{Conjecture}

\newtheorem{quest*}{Question}

\theoremstyle{definition}
\newtheorem{defn}[thm]{Definition}
\newtheorem{defns}[thm]{Definitions}

\theoremstyle{remark}

\newtheorem{rem}[thm]{Remark}

%\newtheorem{sch}[thm]{Scholium}
%\makeatletter
%------------------------The bpp stuff
\newtheoremstyle{montobbo}
{7pt}%spacing above
{5pt}%spacing below
{}
{}%
{\bfseries}
{}%
{.5em}
{\thmnumber{\@{#1}{}~\@{#2}.}%
	\thmnote{~{\bfseries#3.}}}

\theoremstyle{montobbo}
	\newtheorem{montobo}[thm]{}
	\newcommand{\bmo}{\begin{montobo}}
	\newcommand{\emo}{\end{montobo}}
	\newtheorem{talk}{Lecture}
	\newcommand{\btalk}{\begin{talk}}
		\newcommand{\etalk}{\end{talk}}
%%%%%%%%%%%%%%%%%%%%%%%%%%%%%%%%%%%%%%%%%%%%%%%%%%%%%%%%%%%%%%%%%%%%%%%%%%%%%%%%%%%%%%%%%%%%%%%%%
	
%\newcommand{\talk}[1]{\bfseries{#1}}

%%%%%%%%%%%%%%%%%%%%%%%%%%%%%%%%%%%%%%%%%%%%%%%%%%%%%%%%%%%%%%%%%%%%%%%%%%%%%%%%%%%
\newenvironment{manualtheorem}[1]{%%%%%%%%%%%Theorem in the introduction without numbers
	\renewcommand\thethm{#1}%https://tex.stackexchange.com/questions/391443/new-theorem-environment-with-manual-theorem-number
	\thm
}{\endthm}%name as you like
\newcommand{\but}[1]{\begin{manualtheorem}{#1}}
\newcommand{\eut}{\end{manualtheorem}}
\newenvironment{manualproposition}[1]{%%%%%%%%%%%Theorem in the introduction without numbers
	\renewcommand\thethm{#1}%https://tex.stackexchange.com/questions/391443/new-theorem-environment-with-manual-theorem-number
	\prop
}{\endprop}%name as you like
\newcommand{\bup}[1]{\begin{manualproposition}{#1}}
	\newcommand{\eup}{\end{manualproposition}}
%\let\c@equation\c@thm
%\def\els@aparagraph[#1]#2{\elsparagraph[#1]{#2\@addpunct{.}}}
%\def\els@bparagraph#1{\elsparagraph*{#1\@addpunct{.}}}
%\makeatother
\numberwithin{equation}{thm}
%\numberwithin{equation}{section}
%\numberwithin{equation}{subsection}
%\renewcommand{\theequation}{\thesection\Alph{equation}}
%\renewcommand{\theequation}{\arabic{equation}}
%--------Environments define-------------------------------------------
\newenvironment{teqn}{\begin{equation}\begin{tikzcd}\displaystyle}{\end{tikzcd}\end{equation}}
\newenvironment{teqn*}{\begin{equation*}\begin{tikzcd}\displaystyle}{\end{tikzcd}\end{equation*}}

\newenvironment{tpic*}{\begin{equation*}\begin{tikzpicture}}{\end{tikzpicture}\end{equation*}}
%\renewenvironment{proof}{\paragraph{\textit{Proof.}}}{fill}

%--------Symbols define -----------------------------------------------
\DeclareMathAlphabet\matheu{U}{eus}{m}{n}
\DeclareMathAlphabet{\mathmm}{U}{mmambb}{m}{n}
\newcommand{\bb}[1]{\mathbb{#1}}
\newcommand{\fr}[1]{\mathfrak{#1}}
\newcommand{\ca}[1]{\mathcal{#1}}
\newcommand{\scr}[1]{\mathscr{#1}}

%---------Begin Define---------------------------------
\newcommand{\comment}[1]{}
\newcommand{\bun}{\begin{itemize}}
\newcommand{\bn}{\begin{enumerate}}
\newcommand{\bt}{\begin{thm}}
\newcommand{\bl}{\begin{lem}}
\newcommand{\bp}{\begin{prop}}
\newcommand{\bc}{\begin{cor}}
\newcommand{\bd}{\begin{teqn}}
\newcommand{\bud}{\begin{teqn*}}	
\newcommand{\bs}{\begin{proof}}
\newcommand{\br}{\begin{rem}}
\newcommand{\bdf}{\begin{defn}}
\newcommand{\bcj}{\begin{conj}}

\newcommand{\eun}{\end{itemize}}
\newcommand{\en}{\end{enumerate}}
\newcommand{\et}{\end{thm}}
\newcommand{\el}{\end{lem}}
\newcommand{\ep}{\end{prop}}
\newcommand{\ec}{\end{cor}}
\newcommand{\ed}{\end{teqn}}
\newcommand{\eud}{\end{teqn*}}
\newcommand{\es}{\end{proof}}
\newcommand{\er}{\end{rem}}
\newcommand{\edf}{\end{defn}}
\newcommand{\ecj}{\end{conj}}

%               -------  Calculus   --------------

%           ------- Linear Algebra ---------------

\DeclareMathOperator{\GL}{GL}

% ---------------Algebra

\DeclareMathOperator{\trdeg}{tr.deg}

\DeclareMathOperator{\codim}{codim}

%           -------- Geometry ---------------
\DeclareMathOperator{\spec}{Spec}

%          ---------- Arithmetic -------------

%\ReDeclareMathOperator{\et}{\Acute{e}t}

%\DeclareMathOperator{\Qproots}{\mathbb{Q}_\text{p}\left(\text{p}^{1/\text{p}^{\infty}}\right)}
%\DeclareMathOperator{\eq}{\xlongrightarrow{\text{eq.}}}
% ---------------------Topology ----------------------------

%---------------------Diagrams--------------------------------

\newcommand{\iso}{\xrightarrow{\sim}}
\newcommand{\isom}{\arrow[r, "\sim"]}
%----------------------Commands define--------------------------------

\makeatletter\makeatother